\documentclass[11pt]{article}
\usepackage{amssymb}
\usepackage{latexsym}
\newtheorem{thm}{Theorem}[section]
\newtheorem{lem}[thm]{Lemma}
\newtheorem{cor}[thm]{Corollary}
\newtheorem{pro}[thm]{Proposition}

\newtheorem{defi}[thm]{Definition}

\newcommand{\gm }{\Gamma }
\newcommand{\lon }{\longrightarrow }
\newcommand{\be }{\begin{eqnarray*}}
\newcommand{\ee }{\end{eqnarray*}}

\newcommand{\per }{\backl }

\input amssym.def
\input amssym

\newcommand{\pf}{\noindent{\bf Proof.}\ }
\newcommand{\qed}{\begin{flushright} $\Box$\ \ \ \ \ \
                  \end{flushright}}
\newcommand{\complex}{{\mathbb C}}
\newcommand{\reals}{{\mathbb R}}

\newcommand{\frakg}{{\mathfrak g}}

\newcommand{\backl}{\mathbin{\vrule width1.5ex height.4pt\vrule height1.5ex}}

\newcommand{\calc}{{\cal C}}

\newcommand{\call}{{\cal L}}

\newcommand{\calo}{{\cal O}}

\newcommand{\calx}{{\cal X}}

\newcommand{\smalcirc}{\mbox{\tiny{$\circ $}}}

\def\description label#1{\hfil\bf[#1]\hfil}
\newcommand{\tfrakg}{\frakg \times \frakg \times\frakg}
\newcommand{\omeb}{\omega^{b}}
\newcommand{\difft}{\frac{d}{dt}}
\newcommand{\SSS}{{\cal F}}
\newcommand{\TT}{{\cal T}}
\newcommand{\peso}{pseudo}
\newcommand{\perr}{\per}
\newcommand{\kahler}{k\"ahler }

\def\sdp{\mathbin{\hbox{$\mapstochar\kern-.3333em\times$}}}
\def\pds{\mathbin{\hbox{$\times\kern-.55em\mapstochar\,$}}}

\newcommand{\wed}{\mathbin{\lower1.5pt\hbox{$\scriptstyle{\wedge}$}}}

\let\Tilde=\widetilde

\def\chigh{{\raise1.5pt\hbox{$\chi$}}}
\let\phi=\varphi

\def\til0{\Tilde{0}}

\def\dminus{\raise2pt\hbox{\vrule height1pt width 2ex}\hskip3pt}

\def\pback#1{\mathbin{{{\lower1.2ex\hbox{$\times$}}\atop #1}}}

\def\vlra{\hbox{$\,-\!\!\!-\!\!\!-\!\!\!-\!\!\!-\!\!\!
-\!\!\!-\!\!\!-\!\!\!-\!\!\!-\!\!\!\longrightarrow\,$}}

\def\gpd{\,\lower1pt\hbox{$\longrightarrow$}\hskip-.24in\raise2pt
             \hbox{$\longrightarrow$}\,}

\def\lgpd{\,\lower1pt\hbox{$\vlra$}\hskip-1.02in\raise2pt\hbox{$\vlra$}\,}

\def\llgpd{\,\lower1pt\hbox{$\vvlra$}\hskip-1.3in\raise2pt\hbox{$\vvlra$}\,}

\begin{document}

\title{{\bf Structure  ``Hyper-Lie Poisson" }
}

\author{ PING XU \thanks {Research partially supported by NSF
	grants DMS92-03398 and DMS95-04913.}\\
	  Department of Mathematics\\
		    Pennsylvania State University\\
		    University Park, PA 16802, USA \\
		     {\sf email: ping@math.psu.edu}}

\date{}

\maketitle
\begin{abstract}
The main purpose of the paper is to study hyper\kahler structures
from the viewpoint of symplectic geometry.
We introduce a  notion of hypersymplectic structures  which 
encompasses that of hyper\kahler structures.
Motivated by the work of Kronheimer on (co)adjoint orbits
of semi-simple Lie algebras \cite{Kronheimer:LMS} \cite{Kronheimer:JDG1990},
we define  hyper-Lie Poisson structures associated with
a compact semi-simple Lie algebra and give criterion which implies their existence.
 We study an explicit    example of a hyper-Lie Poisson structure,
in which the moduli spaces of solutions to Nahm's equations
assocaited to Lie algebra $\frak{su}(2)$ are
realized  as hypersymplectic leaves and are related to  the  (co)adjoint orbits
of $\frak{sl}(2, \complex )$.
\end{abstract}

\section{Introduction}

Due to  its  rich structure and close connection
  with gauge theory, hyper\kahler manifolds have
 attracted  increasing interest
 \cite{Atiyah} \cite{AtiyahH:1988} \cite{HKLR}. Roughly speaking,
a hyper\kahler manifold is a Riemannian manifold with   three
compatible complex structures  $I, J$ and $K$. The compatibility
means that $I, J, K$ satisfy the quaternion identities
$I^2 =J^2 =K^2 =IJK =-1$, and the metric is k\"ahlerian
with respect to   $I, J$ and $K$.
While it is  easy to find examples
of K\"ahler manifolds, hyper\kahler manifolds are
in general more  difficult to construct.
The two main  often used routes   are twistor theory \cite{Penrose}
and hyper\kahler reduction  \cite{HKLR}, a   generalization
of Marsden-Weinstein reduction \cite{MarsdenW} in the hyper-context.

To  any    K\"ahler manifold there associates a    symplectic
structure, namely its K\"ahler form.  For a  hyper\kahler manifold, there are  
three  symplectic structures (or equivalently,  Poisson structures)
 compatible  with one another in a certain sense (described in Section 2).  However,
there is  an essential difference 
 between K\"ahler and hyper\kahler manifolds.
 For K\"ahler manifolds,  there might be  more than  one
metric compatible with the same symplectic structure.
In  contrast, hyper\kahler structures are more rigid.  Namely,
 the three symplectic structures completely
 determine the hyper\kahler metric,  whence  the  corresponding
complex structures. 
This   observation   suggests
another   possible way of    constructing   a  hyper\kahler manifold, namely,
by  constructing three compatible symplectic structures
on the manifold.

By focusing on symplectic structures rather than the metrics,
we arrive at a new way to define the hyper\kahler condition.
This leads to our definition of  hypersymplectic manifolds.
More precisely, a hypersymplectic manifold
is a manifold  which  admits    three
symplectic structures 
satisfying the same compatibility condition as  usual hyper\kahler
manifolds.
Hypersymplectic manifolds and their basic properties
will be studied in detail  in Section 2.


A powerful  method of constructing  symplectic manifolds  in symplectic
geometry is by  means of Poisson  manifolds.  Every Poisson
manifold admits  a natural  foliation, called the  symplectic
foliation, whose leaves are all   symplectic  \cite{Weinstein:1983}.
For example, symplectic structures on
 coadjoint orbits  are induced from the Lie-Poisson
structure on the  Lie algebra dual  $\frakg^*$ \cite{Weinstein:1983}.
It is   reasonable to expect that there exist 
  some examples  of
``hyper-Poisson manifolds". By a hyper-Poisson manifold,  we mean a
manifold  which admits three  Poisson structures
 satisfying certain  compatibility condition
 such that each leaf is hypersymplectic.
Instead of developing general theory of hyper-Poisson structures,
 the present paper   will be focused on finding some
interesting   examples. In particular, we will consider
 hyper-Lie Poisson structures, an analogue of  
 Lie-Poisson structures in the hyper-context,
which  are presumably the
simplest and most interesting  examples.

There are at least two reasons that hyper-Lie Poisson structures
are important. The  first one comes from our attempt to understand
 general hyper\kahler reduction.
According to the author's knowledge, a successful reduction
theory only exists at $0$ so far. While  Lie-Poisson
structures (or their symplectic leaves) play
an essential role in general symplectic reduction,
carrying out general hyper\kahler reduction will challenge  our
knowledge of hyper-Lie Poisson structures.  In this aspect,
an open question is

\begin{quote}
{\bf Open Problem: } Suppose that $X$ is a hyper\kahler manifold, on which the Lie group $G$
acts preserving the hyper\kahler structure.  Assume that
the action is ``good"  enough  so that the quotient 
$X/G$ is a manifold. Then the three Poisson structures on $X$
corresponding to the three  K\"ahler forms will descend to
three Poisson structures on  the quotient $X/G$. 
What are  the  relation between these  reduced Poisson structures, 
hyper-Poisson structures
and general hyper\kahler reduction?
\end{quote}

Understanding Kronheimer's recent work on adjoint
orbit hyper\kahler structures provides anther strong  motivation. 
Using gauge theory and infinite
dimensional hyper\kahler reduction,
Kronheimer proved that
certain adjoint orbits of complex semi-simple
Lie algebras admit hyper\kahler metrics   \cite{Kronheimer:LMS} \cite{Kronheimer:JDG1990}.
Later on, his  result was generalized by  Biquard and Kovalev
to  arbitrary  adjoint
orbits \cite{Biquard} \cite{Kovalev}, and has been  used to 
understand  the Kostant-Sekiguchi correspondence
\cite{Vergne}. However,
the hyper\kahler metrics and symplectic structures
obtained in this way are quite  mysterious and   elusive.
Recall that an intrinsic reason for  each coadjoint orbit  to admit
a symplectic structure is that the  dual of any Lie algebra is a Poisson manifold,
and each coadjoint orbit happens  to be a symplectic leaf of this  Lie-Poisson
structure.
Inspired by this fact as well as the results of Kronheimer and others,
 it is quite reasonable
 to expect that there exists a hyper-Lie Poisson structure
such that the  orbits studied by Kronheimer et al   occur
as  its hyper-symplectic leaves.  Then, this will provide  us a natural source and
symplectic explanation
 for those hyper\kahler structures on adjoint orbits.
To explore the connection between   the work of  Kronheimer
  and   symplectic
geometry was indeed     the initial  motivation for us to consider hyper-Lie
Poisson structures. 

In this paper, as an example, we will consider in detail a hyper-Lie
Poisson structure associated with $\frakg =\frak{su}(2)$. 
In the meantime, we will take some tentative steps toward hyper-Lie
Poisson structures  associated with general  compact semi-simple
Lie algebras. For this purpose, we will 
keep the  discussion  general  from Section 2 through 
Section 3 while the last two sections will be  devoted to the 
special case  $\frakg =\frak{su}(2)$.

To  explain our approach, we need to rephrase the
definition of hyper\kahler manifolds in a  way slightly different from the
literature.
Note that there is in fact  no preferred choice of complex
structures on a hyper\kahler manifold.
The  bundle maps
$I', J'$ and $K'$ given by $(I' , J', K')=
(I, J, K)O$,  for any orthogonal matrix $O\in SO(3)$,
will satisfy exactly the same  quaternion
 relations.  Therefore the map:  $O\in SO(3)  \lon I'$ assigns a
complex structure to every  orthogonal matrix in $SO(3)$.
In particular, under such an assignment,
$I$, $J$ and $K$ are the complex structures corresponding to 
the identity matrix,
the   matrix of  the 
cyclic permutation: $\{ e_{1}, e_{2}, e_{3}\}\lon \{ e_{2}, e_{3}, e_{1}\}$ and the  matrix of 
the  cyclic permutation:
$\{ e_{1}, e_{2}, e_{3}\}\lon\{e_{3}, e_{1}, e_{2}\}$, respectively.
Since a matrix in  $SO(3)$  can  be naturally   identified  with
a standard orthonormal  basis  in $\frak{su}(2)$  (such  are also
called frames in this paper),
intrinsically  we can think of  a hyper\kahler
manifold as a manifold with   a family of complex structures (or equivalently
symplectic structures),
 parameterized by frames. This point of view is different
from the conventional one, in which  complex structures (or symplectic
structures)  on a hyper\kahler manifold are  considered to be
parameterized by the unit sphere  $S^2$.  This is   the crux
in our approach.

Now  the question arises  in  which space  
 a hyper-Lie Poisson structure should  live. To answer this question, 
we first recall that  a Lie-Poisson space $\frakg^*$
emerges as the  target space of a   momentum mapping  of a symplectic
homogeneous $G$-space.
A momentum mapping of a homogeneous hyper\kahler $G$-space $X$ is usually
considered, when the  three complex structures $I, J, K$
are chosen, as 
 a map from  $X$ to $ \frakg^* \times \frakg^* \times
\frakg^*$   \cite{HKLR}.
However, when  $I, J, K$
are replaced by any other three complex structures $I', J', K'$
related by an orthogonal matrix in $SO(3)$ as
described earlier, the corresponding
momentum mapping  changes accordingly.
Intrinsically, the momentum mapping of a homogeneous
hyper\kahler manifold
should  therefore be considered as a  map $X\lon L(\frakg, \frak{su}(2))$,
where $L(\frakg,   \frak{su}(2) )$  is 	the space of all linear maps from
$\frakg$ to $ \frak{su}(2)$.
It is therefore  reasonable to
    expect that  $L(\frakg, \frak{su}(2) )$, 
    as the target space of  the momentum mapping of a homogeneous
hyper\kahler manifold,   should carry  
a  hyper-Lie Poisson  structure.
 Another possibly  useful way to    think  of    $L(\frakg, \frak{su}(2))$
as a natural generalization of $\frakg^*$
 is to  note that   this space  is obtained 
 from $ \frakg^* = L(\frakg, \reals )$
 by replacing $\reals$ by $\frak{su}(2)$.

An elegant way of obtaining a family of Poisson structures
on the space  $L(\frakg,   \frak{su}(2) )$
goes as follows. The space  $L(\frakg,   \frak{su}(2) )$ 
can be identified with $\frakg^* \times \frakg^* \times
\frakg^*$, once  an orthonormal basis of $ \frak{su}(2)$, i.e., a frame,
 is fixed.
If we  can
 define a Poisson structure  $\pi$ on the space  $\frakg^* \times \frakg^* \times
\frakg^*$, by pulling back $\pi$ to $ L(\frakg,   \frak{su}(2) )$ under such
an identification, we then  obtain a  Poisson structure on the space 
$ L(\frakg,   \frak{su}(2) )$.
 This construction  in fact enables us to
obtain a family of Poisson structures   simply by  varying  the frames.
  Throughout the paper, we shall identify $\frakg^*$ with $\frakg$ via  the Killing form, hence $\frakg^* \times \frakg^* \times
\frakg^*$ with $\tfrakg$ for simplicity. We note that 
any bivector field   on $\tfrakg$ is determined by
its corresponding  brackets of  all linear functions 
$l_{\xi}^{i}$, where $l_{\xi}^1 (a, b, c)=<\xi , a>$ etc.
Therefore, to define  the  Poisson tensor $\pi$, it suffices
to find its corresponding brackets of linear
functions. There are 
 certain natural   conditions that  the Poisson tensor $\pi$ has to satisfy.
One of them is   that both projections $pr_{12}, pr_{13}:
\tfrakg \lon \frakg^{\complex} $ be  Poisson maps.
 Here $pr_{12}, pr_{13}$
 are the maps defined by $pr_{12}(a, b, c)=a+ib$
 and  $pr_{13}(a, b, c) =a+ic$ respectively, and meanwhile
 $\frakg^{\complex} $
 is identified with its dual and  equipped with the Lie Poisson structure
 as a real Lie algebra. By using this condition for $\pi$, 
one can easily write down its  brackets of all linear functions
except for  the bracket   of $l_{\xi}^2 $ and $\l_{\eta}^3$,
which should  correspond to a $S^2 (\frakg )$-valued
function $A$ on $\tfrakg$. Therefore, the entire
problem reduces to that of finding a suitable 
function $A$.

Section 3  is devoted to  a detailed discussion of this  problem as
well as an investigation of  when $\pi$
defines a Poisson tensor. Moreover, we will derive
the  criterion  on $A$  that
 the induced family of Poisson structures on   $L(\frakg,   \frak{su}(2) )$ is 
compatible and have the desired  properties as outlined above.

Section 4 is devoted to the study of the case $\frakg=\frak{su}(2)$,
where a satisfactory function $A$ is explicitly constructed
on an open submanifold  of  $L(\frakg,   \frak{su}(2) )$.
The induced hypersymplectic
foliation is also explicitly  described,  and  a complete
set of casimirs is obtained. The corresponding  \peso-metric on each leaf
is also computed,  and in fact it is shown that all 
hypersymplectic leaves are hyper\kahler.

Section 5 is a continuation of Section 4, where the  hyper-Poisson
structure is   extended to a certain critical set $C$.
It is shown that hyper-symplectic leaves of this extended
hyper-Lie Poisson structure are  diffeomorphic
to 
(co)adjoint orbits of $\frak{sl}(2, \complex )$.
In this way, we obtain a symplectic proof for the existence
of hyper\kahler structures on (co)adjoint orbits of  $\frak{sl}(2, \complex )$.
Although $\frak{sl}(2, \complex )$ is the simplest semi-simple Lie
algebra, the existence of hyper\kahler structures on its
(co)adjoint orbits is somewhat already nontrivial (see \cite{Hitchin}).




{\bf Acknowledgments.}  
The author  would like to thank Jorgen  Andersen,  Sean Bates, Olivier
Biquard,  Ranee Brylinski,
Hansj\"{o}rg Geiges, Victor Guillemin, 
Nigel  Hitchin,  Peter  Kronheimer,  Oscar Garcia-Prada, Yvette Kosmann-Schwarzbach
 and Alan Weinstein for useful discussions and email correspondences.
Special thanks go to Ludmil Katzarkov and Tony Pantev
for calling  his attention to the papers  \cite{Kronheimer:LMS}
\cite{Kronheimer:JDG1990}.
Most of this work was carried out when the author was a member of MSRI, where
the research was supported by NSF Grant DMS90-22140.
In addition, the author wishes to thank
the  Center of Mathematics at  Zhejiang University, China, 
for its hospitality while part of this work being done.
  
\section{Hypersymplectic structures}
The purpose of the present  section is to introduce a  notion called
hypersymplectic structures,
 which includes hyper\kahler  manifolds as a  special case
and  is much  more natural  from the viewpoint of
 symplectic geometry.
Our definition of hyper\kahler structures here  is slightly  different from  the
one in the literature, where complex structures and metrics  have received
much more attention. Our   interests in this paper  mainly  lie
in symplectic forms and their Poisson tensors.

By  $\Omega^2_{+} (S)$, we    denote the space of nondegenerate
 2-forms on a manifold $S$, and $\gm_{+}(\wedge^2 TS)$
the space of non-degenerate bivector fields on $S$. By
$\kappa$, we denote the  map $\Omega^2_{+} (S)\lon \gm_{+}(\wedge^2 TS)$,
 which is the inversion when  elements in $\Omega^2_{+} (S)$ and
$\gm_{+}(\wedge^2 TS)$ are   considered as bundle maps. For any
$\omega\in \Omega^2_{+} (S)$, we write $\omega^{-1}=\kappa (\omega )\in
\gm_{+}(\wedge^2 TS)$.
A $\frak{su}(2)$-valued 2-form $\Omega$ on  $S$ is said to be
nondegenerate if the form $\omega_{\xi}=<\xi, \Omega>$
is nondegenerate for any non-trivial $\xi \in \frak{su}(2)$,
where the pairing is with respect to the Killing form on $\frak{su}(2)$.
\begin{defi}
A hypersymplectic structure on a manifold $S$ is a closed non-degenerate
$\frak{su}(2)$-valued two-form $\Omega$ such that $\xi \lon \omega_{\xi}^{-1} $
 is   linear when $ \|\xi \|=1$.
 That is, for any unit vectors $\xi, \xi_{1}, \xi_{2}\in \frak{su}(2)$
satisfying  $\xi =k_{1}\xi_{1}+k_{2}\xi_{2} $,  
\begin{equation}
\omega_{\xi}^{-1}= k_{1}\omega_{\xi_{1}}^{-1}+k_{2} \omega_{\xi_{2}}^{-1} .
 \end{equation}
\end{defi}

Suppose that $\Omega\in \Omega^2 (S)\otimes \frak{su}(2)$ is a hypersymplectic
structure on $S$. By choosing an orthonormal basis $\{e_{1}, e_{2}, e_{3}\}$
of $\frak{su}(2)$,   also called {\em a frame} in
the sequel, $\Omega$ can be written as 
$\Omega=\omega_{1}e_{1}+\omega_{2}e_{2}+\omega_{3}e_{3}$, where
$\omega_{1}, \omega_{2}$ and $\omega_{3}$ are  symplectic forms on $S$.
By $\pi_{1}, \pi_{2}$ and $\pi_{3}$ we denote
their corresponding bivector fields on $S$.

As usual, for each $i$,  $\omega_{i}^{b}$ denotes
the bundle map $TS\lon T^* S$ given by $<\omega_{i}^{b}(v), u>
=\omega_{i}(v, u)$,  $\forall u, v\in TS$,  $X^{i}_{f}$
the vector field on $S$ defined  by $\omega_{i}^{b} (df)$ for any
 $f\in C^{\infty}(S)$,
and $\{f, g\}_{i}=X^{i}_{f} g$ for any $f, g\in C^{\infty}(S)$.


\begin{pro}
Suppose that $\Omega=\omega_{1}e_{1}+\omega_{2}e_{2}+\omega_{3}e_{3}
\in \Omega^2 (S)\otimes \frak{su}(2)$ is a hypersymplectic structure
on $S$. Then,

(1).
\begin{equation}
\label{eq:compatible}
[\omega_{i}^b\smalcirc (\omega_{j}^{b})^{-1}]^2=-1, \ \ \ \ \forall i\neq j,
\end{equation}
where $1$ denotes the identity  map  $TS\lon TS$.

(2). $$[\pi_{i}, \pi_{j}]=0, \mbox{  for any} \ i, j, $$
where the bracket is the Schouten bracket on multivector fields \cite{K-SM:1990}.

\end{pro}
\pf For any $k_{1}, k_{2}, k_{3}$ such that $k_{1}^2+k_{2}^2
+k_{3}^2 =1$, $k_{1}e_{1}+k_{2}e_{2}+ k_{3}e_{3}$ is a
unit vector in $\frak{su}(2)$. Thus,  it follows from definition that
$$(k_{1}\omega_{1}^b + k_{2}\omega_{2}^b  +k_{3}\omega_{3}^b )^{-1}=
k_{1} (\omega_{1}^b )^{-1}+k_{2}(\omega_{2}^b )^{-1}+k_{3}(\omega_{3}^b )^{-1}.$$
That is,
\be
1&=&(k_{1}\omega_{1}^b + k_{2}\omega_{2}^b   +k_{3}\omega_{3}^b )
[k_{1} (\omega_{1}^b )^{-1}+k_{2}(\omega_{2}^b )^{-1}+k_{3}(\omega_{3}^b )^{-1}]\\
&=&(k_{1}^2+k_{2}^2 +k_{3}^2)+k_{1}k_{2}[\omega_{1}^b  (\omega_{2}^b )^{-1}
+\omega_{2}^b (\omega_{1}^b )^{-1}] \\
&&\ \ \ \ 
 +k_{2}k_{3}[\omega_{2}^b (\omega_{3}^b )^{-1}
 +\omega_{3}^b (\omega_{2}^b )^{-1}]
 +k_{1}k_{3}[\omega_{1}^b (\omega_{3}^b )^{-1}
 +\omega_{3}^{b} (\omega_{1}^b )^{-1}].
 \ee
 Equations (\ref{eq:compatible}) thus follow immediately.

 Also, from the argument above, we see that $k_{1}\pi_{1}+ k_{2}\pi_{2}
 +k_{3} \pi_{3}$ is still a Poisson tensor for any $(k_{1}, k_{2}, k_{3})$
 in the 2-sphere, hence for arbitrary $k_{1}, k_{2}, k_{3}$ as well.
 It thus follows that $[\pi_{1}, \pi_{2}]=[\pi_{2}, \pi_{3}]
 =[\pi_{1}, \pi_{3}]=0$. \qed

An immediate consequence is the following:
 \begin{cor}
 Let $S$ be  a hypersymplectic manifold with hypersymplectic
 form $\Omega \in \Omega^{2}(S)\otimes \frak{su}(2)$. Then for
 any $\xi ,  \eta \in \frakg$,
 $$[\pi_{\xi}, \pi_{\eta}]=0, $$
 where $\pi_{\xi}=\omega_{\xi }^{-1}$ and $\pi_{\eta}=\omega_{\eta }^{-1}$.
\end{cor}

Sometimes the following equivalent version
 is more often used.

\begin{pro}

Equations (\ref{eq:compatible}) are equivalent to
\begin{equation}
\label{eq:compatible1}
\omega_{i}(X_{f}^j , X_{g}^j )=\{f, g\}_{i}, \ \ i\neq j, \mbox{ for any }
f, g\in C^{\infty} (S).
\end{equation}
\end{pro}
\pf It is quite obvious that
\be
[\omega_{i}^b\smalcirc (\omega_{j}^{b})^{-1}]^2=-1 &\Longleftrightarrow &
(\omega_{j}^{b})^{-1} \omega_{i}^b (\omega_{j}^{b})^{-1}= -(\omega_{i}^{b})^{-1}\\
&\Longleftrightarrow&
<(\omega_{j}^{b})^{-1} \omega_{i}^b (\omega_{j}^{b})^{-1}df, dg>
=-<(\omega_{i}^{b})^{-1} df, dg>,
\ee
 which is equivalent to
$$\omega_{i}(X_{f}^j , X_{g}^j )=\{f, g\}_{i}. $$ \qed

In fact, Equations (\ref{eq:compatible}),  or equivalently Equations (\ref{eq:compatible1}),
 are also  sufficient to construct a
 hypersymplectic structure on $S$.

 \begin{pro}
 Suppose that  $S$  is a manifold with three symplectic structures
 $\omega_{1}, \omega_{2}$ and $\omega_{3}$ such that
  Equations (\ref{eq:compatible}) (or Equations (\ref{eq:compatible1}))
  hold. Then for any orthonormal basis $\{e_{1}, e_{2}, e_{3}\}$
of $\frak{su}(2)$,  the $\frak{su}(2)$-valued 2-form
$\Omega=\omega_{1}e_{1}+\omega_{2}e_{2}+\omega_{2}e_{2}$
  defines a hypersymplectic structure on $S$.
  \end{pro}

The proof  is quite straightforward,  and  is
left  for the reader.

  To each  hypersymplectic manifold, we  associate
  a   natural \peso-metric as we will see below.

Let  $g: TS\lon T^*S$ be the bundle map
given by 
\begin{equation} 
\label{eq:g}
g=\omega_{3}^{b}(\omega_{1}^{b})^{-1}\omega_{2}^{b}, 
\end{equation}
and $I$, $J$, $K$ the bundle maps form $TS$ to itself defined by
\begin{equation}
I=g^{-1}\omega_{1}^{b},\ \  J=g^{-1}\omega_{2}^{b}\ \  \mbox{ and }
K=g^{-1}\omega_{3}^{b}.
\end{equation}

\begin{thm}
\begin{enumerate}
\item $g$ is a \peso-metric on $S$; 
\item $g$ can  also be  written as $g=\omega_{1}^{b}(\omega_{2}^{b})^{-1}\omega_{3}^{b}$,
or  $g=\omega_{2}^{b}(\omega_{3}^{b})^{-1}\omega_{1}^{b}$;
\item  $I, J, K$ satisfy the quaternion
relation:
$$I^2=J^2=K^2=IJK=-1.$$
\end{enumerate}
\end{thm}
\pf (1).  It follows from Equations  (\ref{eq:compatible}), by taking the dual,
that
\begin{equation}
\label{eq:compatible2}
[(\omega^b_{j})^{-1}  \omega^b_{i} ]^2=-1, \ \ \forall i\neq j.
\end{equation}
Hence, by  using Equations (\ref{eq:compatible})  
and (\ref{eq:compatible2}) repeatedly, we have 
$g^*=-\omega_{2}^{b}(\omega_{1}^{b})^{-1}\omega_{3}^{b}
=[\omega_{1}^{b}(\omega_{2}^{b})^{-1}]\omega_{3}^{b}
=-\omega_{1}^{b}[(\omega_{3}^{b})^{-1}\omega_{2}^{b}]
=[\omega_{3}^{b}(\omega_{1}^{b})^{-1}]\omega_{2}^{b}
=g$.
 That is, $g$ is symmetric. Furthermore, it is
 evident  that $g$ is nondegenerate since 
$\omega_{1}, \omega_{2}, \omega_{3}$ are all nondegenerate.

(2). It follows from Part  (1) that $g=g^*=-\omega_{2}^{b}(\omega_{1}^{b})^{-1}
\omega_{3}^{b}=\omega_{1}^{b}(\omega_{2}^{b})^{-1}\omega _{3}^{b}$.
The other  equation can be obtained in  a similar way.

(3).  Using Part   (2), we have $I=g^{-1}\omega_{1}^{b}=
[\omega_{1}^{b}(\omega_{2}^{b})^{-1}\omega _{3}^{b}]^{-1}\omega_{1}^{b}=
(\omega _{3}^{b})^{-1}\omega_{2}^{b}$.
Similarly,  $J=(\omega _{1}^{b})^{-1}\omega_{3}^{b}$ and
$K=(\omega _{2}^{b})^{-1}\omega_{1}^{b}$. 
Therefore, $I^2=J^2=K^2 =-1$.  Furthermore,
\be
IJ&=& (\omega _{3}^{b})^{-1}\omega_{2}^{b} (\omega _{1}^{b})^{-1}
\omega_{3}^{b}\\
&=&[\omega_{1}^{b}(\omega_{2}^{b})^{-1}\omega _{3}^{b}]^{-1}\omega_{3}^{b}\\
&=&g^{-1}\omega_{3}^{b}\\
&=&K.
\ee
This concludes  the proof. \qed 

It looks as if  our definition of $g$    depends on a particular
choice of the  frame. However, the following theorem indicates
that 
$g$  is in fact independent of frames up to a sign.  

\begin{thm}
If two frames are of the same orientation, then their corresponding
 pseudo-metrics  coincide.
\end{thm}
\pf  Let $\SSS=\{e_{1}, e_{2}, e_{3}\}$ and $\TT=\{e_{1}', e_{2}', e_{3}'\}$
be any two frames.  Suppose that $(e_{1}', e_{2}', e_{3}' )=(e_{1}, e_{2}, e_{3})O$ for some orthogonal matrix $O\in SO(3)$.
Let $I, J, K$  be the induced almost  complex structures on $S$  corresponding
to the frame $\SSS$. We define $I', J',  K'$ by
the equation:
$$(I', J', K')=(I, J, K)O. $$ 
It follows from the    quaternion relation of $I, J, K$
that $I', J',  K'$ also satisfy the  same relation.
Since $ \omega_{1}e_{1}+\omega_{2}e_{2}+\omega_{2}e_{2}=
\omega_{1}'e_{1}'+\omega_{2}'e_{2}'+\omega_{2}'e_{2}'$,
it follows that $(\omega_{1}', \omega_{2}', \omega_{3}')=
(\omega_{1}, \omega_{2},  \omega_{3})O$. Hence, we have
\begin{equation}
I'=g^{-1}(\omega_{1}')^{b}, J'=g^{-1}(\omega_{2}')^{b} \mbox{ and }
K'=g^{-1}(\omega_{3}')^{b}.
\end{equation}
By using the quaternion  relation of $I', J',  K'$, we can
easily deduce that $g= 
(\omega_{3}')^{b}((\omega_{1}')^{b})^{-1}(\omega_{2}')^{b}$,
which is $g'$ by definition. \qed

Because of this result, we call $g$  the \peso-metric
associated to the hypersymplectic structure despite of 
an ambiguity of  signs. In particular, if $g$ is 
positive (or negative) definite, the hypersymplectic structure becomes
 hyper\kahler.
 We refer the reader
to \cite{Atiyah}  \cite{AtiyahH:1988} \cite{Hitchin}
 \cite{HKLR} for the background on the subject of 
  hyper\kahler structures. 

The following result is well-known for hyper\kahler structures,
and  is however still valid in our general context.
Readers can find a proof in, for example, \cite{AtiyahH:1988}.
For completeness, we outline a proof here.
\begin{thm}
If $ S$ is a hypersymplectic manifold, then   the almost complex
structures $I, J, K$ corresponding to any frame  are integrable.
\end{thm}
\pf 
For any vector fields $X, Y\in \calx (S)$,
$$\omega_{2}(X, Y)=g(JX, Y)=g(KIX, Y)=\omega_{3}(IX, Y).$$
 Hence, we have the relation:
\begin{equation}
\label{eq:omega23}
X\per \omega_{2}=IX \per \omega_{3}. 
\end{equation}

It follows that a complex vector field $X$ is of type $(1, 0)$
with respect to $I$ iff 
\begin{equation}
\label{eq:omega23i}
X\per \omega_{2}=iX \per \omega_{3}.
\end{equation}

Suppose that $X, Y$ are complex vector fields of type
$(1, 0)$. In order to show that  $I$ is integrable
it suffices to show that their bracket 
$[X, Y]$ is of type $(1, 0)$ according to
the Newlander-Nirenberg theorem.
However,
$$[X, Y]\per \omega_{2}=L_{X}(Y\per \omega_{2})-Y\per ( L_{X}\omega_{2}). $$
Now,
\be
 L_{X}\omega_{2}&=&(d\iota_{X}+\iota_{X}d)\omega_{2}\\
 &=&d(i\iota_{X}\omega_{3})\\
  &=&iL_{X}\omega_{3}, 
  \ee
  and from Equation (\ref{eq:omega23i}),
  $$Y\per \omega_{2}=i(\iota_{Y}\omega_{3}). $$
  Thus,

  $$[X, Y]\per \omega_{2}=i[L_{X}(\iota_{Y}\omega_{3})-\iota_{Y}(L_{X}\omega_{3})]
=i([X, Y]\per \omega_{3}). $$

  Thus, $I$ is integrable. Similarly, $J$ and $K$ are also
  integrable. \qed
{\bf Remark.}\ \ It is not difficult to check that $(\pi_{2}, I)$, 
$(\pi_{3}, I)$ and all other similar pairs are Poisson-Nijenhuis structures
in the sense  of Kosmann-Scharzbach and Magri \cite{K-SM:1990}.
Poisson-Nijenhuis structures are introduced by
 Kosmann-Scharzbach and Magri in the study of integrable
systems. Therefore, it would be very interesting to explore
the relation between hypersymplectic manifolds and integrable
systems.  

\section{Hyper-Lie Poisson structures}
This section is devoted to the introduction of
hyper-Lie Poisson structures. The main idea is to  define
on a suitable space $M$  a family of Poisson structures
parameterized by frames, which will coherently 
depend on the parameterization  in a proper sense.
We shall analyses the condition under which   the induced 
symplectic  foliations are independent of frames
so that each leaf becomes hypersymplectic. 

To begin with,
let $\frakg$ be a semisimple Lie algebra with Killing form 
$<\cdot , \cdot >$, and $A$ a  $S^2 (\frakg )$-valued function on
$\tfrakg$,   where
$S^2 (\frakg )$ denotes the space of second order symmetric
tensors on $\frakg $.  We note that any element in $S^2 (\frakg )$
can be   naturally considered  as a symmetric bilinear form on $\frakg$.
So contracting with  $\xi, \eta \in\frakg$,  there corresponds to a
function on $\tfrakg $:  $A_{\xi, \eta}=\xi \per A\perr \eta$.
For any $\xi \in \frakg$, we denote by  $l^{1}_{\xi }, l^{2}_{\xi }, l^{3}_{\xi
}$
the  linear functions on  $\frakg \times \frakg \times \frakg$
defined by $l^{1}_{\xi }(a, b, c)=<\xi , a>$, etc.
Our first step is to define a Poisson structure
on $\tfrakg$. In order to do so
 (or even just to define  a bivector field  on $\tfrakg$),
 it  suffices  to define its corresponding brackets among all 
linear functions $l_{\xi}^{i}, \
i=1, 2, 3$,  since they span the function space $C^{\infty}( \tfrakg ) $.

\begin{defi}
Let $A$ be a $S^2 (\frakg )$-valued function on $\tfrakg$.
The following bracket defines a bivector field $\pi $ on  $\tfrakg$.

\be
\{l^{1}_{\xi}, l_{\eta }^1 \}&=&l^{1}_{[\xi, \eta]}\\
\{l^{1}_{\xi}, l_{\eta }^2 \}&=& \{l^{2}_{\xi}, l_{\eta }^1 \}=
l^{2}_{[\xi, \eta]}\\
\{l^{2}_{\xi}, l_{\eta }^2 \}&=&-l^{1}_{[\xi, \eta]}\\
\{l^{1}_{\xi}, l_{\eta }^3 \}&=&\{l^{3}_{\xi}, l_{\eta }^1 \}=
l^{3}_{[\xi, \eta]}\\
\{l^{3}_{\xi}, l_{\eta }^3 \}&=&-l^{1}_{[\xi, \eta]}\\
\{l^{2}_{\xi}, l_{\eta }^3 \}&=& -\{l^{3}_{\xi}, l_{\eta }^2 \}
=A_{\xi, \eta}.
\ee
\end{defi}

Let $G$ be a compact Lie group with Lie algebra $\frakg$.
Then $G$ acts on $\tfrakg$   diagonally, with adjoint
action on each factor.

 \begin{pro}
 \label{pro:A-equi}
 The following are equivalent:
 \begin{enumerate}
 \item The bivector field $\pi$ is $G$-invariant;
 \item the map  $A: \tfrakg \lon S^2 (\frakg )$ is $G$-equivariant,
 where $G$ acts on $ S^2 (\frakg )$ by the  adjoint action;
 \item for any $\xi, \eta, \zeta \in \frakg$, 
 $$\hat{\xi}A_{\eta, \zeta}=A_{[\xi, \eta ], \zeta}+A_{\eta ,[\xi , \zeta ]}.$$
 \end{enumerate}
 \end{pro}   
 \pf That (1) and (2) are equivalent is quite evident.

 $(2)\Longleftrightarrow (3)$.
 Suppose that   $A$ is $G$-equivariant. That is,  $A(gx)=Ad_{g} A(x)$.
 It follows, by taking derivative,  that for any $\xi \in \frakg$,
 $\hat{\xi}A=ad_{\xi }A$. Hence for any $\eta, \zeta \in \frakg$,
 \be
 \hat{\xi}A_{\eta, \zeta}&=&\eta \per \hat{\xi}A \perr  \zeta\\
 &=&\eta \per ad_{\xi }A \perr  \zeta\\
 &=&A_{[\xi, \eta ], \zeta}+A_{\eta ,[\xi , \zeta ]}.
 \ee
 The converse  is   also  true  by using  the same
argument   backwards. \qed

The following theorem gives a necessary
and  sufficient   condition for the 
 bivector field $\pi$ to be  a Poisson tensor.
As usual, for any $f\in C^{\infty}(\tfrakg )$ we write $X_{f}$  for 
the vector field $\pi^{\#}(df ) $.

\begin{thm}
\label{thm:conditionA}
$\pi$ is a Poisson tensor iff $A$ is equivariant
and at any point $(a, b, c)\in \tfrakg$,

\begin{eqnarray}    
\xi \per X_{l^2_{\eta }}A-\eta  \per X_{l^2_{\xi }}A&=&[c, [\xi ,\eta ]],
\label{eq:conditionA} \\
\xi \per X_{l^3_{\eta }}A-\eta  \per X_{l^3_{\xi }}A&=&[b, [\eta ,\xi ]],
\label{eq:conditionB}
\end{eqnarray}
 for any $\xi , \eta \in   \frakg $.

\end{thm}
\pf  By $pr_{12}$,   we denote the
projection $\tfrakg \lon \frakg^{\complex}$ given by
$pr_{12} (a, b, c)=a+ib$. Similarly,  $pr_{13}$ denotes the  projection
from $\tfrakg$ to $ \frakg^{\complex}$ given by
$pr_{13} (a, b, c)=a+ic$. It is simple to see that
$Tpr_{12} \pi=Tpr_{13} \pi=\pi^{\complex}$,
where $\pi^{\complex }$ is the       Lie-Poisson   tensor
on $\frakg^{\complex}$, which is 
  identified with its dual
 as  a  real Lie algebra.
Hence the Jacobi identity 
$$\{\{f_{1}, f_{2}\}, f_{3}\} +c.p. =0 $$
holds if $f_i, i=1, 2, 3, $ are linear functions of the  form
$l_{\xi}^1, l_{\xi}^2$, or of the  form
 $l_{\xi}^1,  l_{\xi}^3$, where $c.p.$ stands for the
cyclic permutation.  It remains to check the following three cases: 
(1) $f_{1}=l_{\xi}^1, f_{2}=l_{\eta}^2$ and $f_3 =l_{\zeta}^3$;
(2) $f_{1}=l_{\xi}^2, f_{2}=l_{\eta}^2$ and $f_3 =l_{\zeta}^3$; and
(3)$f_{1}=l_{\xi}^3, f_{2}=l_{\eta}^3$ and $f_3 =l_{\zeta}^2$.

It is simple to see that the Jacobi identity in Case (1)
is equivalent to that $A$ is
$G$-equivalent according to 
Proposition \ref{pro:A-equi}.
As for Case (2),
\be
\{\{l_{\xi}^2, l_{\eta}^2\}, l_{\zeta}^3\}+c.p.&=&
-\{l^1_{[\xi, \eta ]}, l_{\zeta}^3\}+\{A_{\eta, \zeta}, l_{\xi}^2\}
+\{-A_{\xi, \zeta }, l_{\eta}^2\}\\
&=&-l^3_{[[\xi, \eta ], \zeta ]}-X_{l_{\xi}^2}A_{\eta, \zeta}
+X_{l_{\eta}^2}A_{\xi, \zeta} \\
&=&-<[c, [\xi ,\eta ]], \zeta>-<\eta \per X_{l_{\xi}^2}A ,\zeta >
+<\xi \per X_{l^2_{\eta }}A, \zeta >.
\ee 
Thus, the Jacobi identity follows  iff
Equation (\ref{eq:conditionA})  holds. Similarly,
Equation (\ref{eq:conditionB}) is equivalent to
the Jacobi identity for Case (3). \qed

In the proof above, we have in fact shown the
following:

\begin{pro}
\label{pro:poisson}
Assume that both  Equation (\ref{eq:conditionA})
and Equation (\ref{eq:conditionB})  hold for any
$\xi, \eta \in \frakg$.
Then, both $pr_{12}: \tfrakg\lon \frakg^{\complex}$ and
$pr_{13}: \tfrakg\lon \frakg^{\complex}$ are Poisson maps,
where $\frakg^{\complex}$, identified with its dual,
is equipped with the Lie Poisson structure 
 as a real Lie algebra.
\end{pro}

Consider the space  $M=L(\frakg , \frak{su} (2))$, which consists of
all linear maps from $\frakg $ to  $\frak{su} (2)$. Here $\frak{su} (2)$
is  only considered  as a  vector space without using any
Lie algebra structure. $M$ admits a natural   $G$-action induced 
from   the adjoint action on $\frakg$.
Whenever  a frame, i.e.,  an orthonormal basis $\{e_{1}, e_{2}, e_{3}\}$
of $\frak{su} (2)$,  is chosen, $M$ is
identified with $\frakg^* \times \frakg^* \times\frakg^*$,
which can  also  be identified with $\frakg \times \frakg \times \frakg$
using the Killing form on $\frakg$. We shall denote such
an identification $M\lon \tfrakg$ by $\Psi_{\SSS}$.
From now on, we will always identify $\frakg$ with its
dual $\frakg^*$. Using the map $\Psi_{\SSS}$, the Poisson 
structure $\pi$  on $\tfrakg$ is pulled back  to a Poisson structure $\pi_{\SSS}$
on $M$. When  the choice of   frames $\SSS$ varies, we  thus obtain a
family of Poisson structures on $M$ parameterized by
frames. This is the very  structure we are
interested in.
Corresponding to any  frame $\SSS=\{e_{1}, e_{2}, e_{3}\}$, there exist
two frames $\SSS_{2}$ and  $\SSS_{3}$ obtained by
the cyclic permutations: 
$ \{e_{2}, e_{3}, e_{1}\}$, and $\{e_{3}, e_{1}, e_{2}\}$, respectively. We also often  use $\SSS_{1}$ to denote $\SSS$.
The Poisson structures corresponding to $\SSS_{1}$,  $\SSS_{2}$ and  $\SSS_{3}$
are denoted by  $\pi_{\SSS_{1}}, \pi_{\SSS_{2}}$ and $\pi_{\SSS_{3}}$,
respectively. 

By choosing a frame, any
 vector-valued function  $F$ on $\tfrakg$ can be
pulled back to a function  $F_{\SSS}$ on $M$ via the
map  $\Psi_{\SSS}$.
If furthermore,  $F$ is invariant under the action of $O(3)$,
 where
$O(3)$   acts on $\tfrakg$  by $(a, b, c)\lon (a, b, c)O$
for any $(a, b, c)$ and $O\in O(3)$, 
$F_{\SSS}$ is independent of the frame  $\SSS$, and therefore can be
considered as a well-defined function on $M$. However, in most
cases,  the  function $F$ is only invariant under the $SO(3)$-action.
In this case, the pull back  $F_{\SSS}$ depends on the orientation
of $\SSS$. Whenever the orientation of frames is fixed, we
shall still get a well-defined function on $M$.   
In the sequel, we  shall always assume that the $S^2 (\frakg )$-valued
function $A$ on $\tfrakg$ is invariant under the $SO(3)$-action,
and therefore can be considered as a function on $M$ when an
orientation of frames is fixed.

 \begin{thm}
\label{thm:global}
 $\Pi= \pi_{\SSS_{1}}e_{1}+\pi_{\SSS_{2}} e_{2}+\pi_{\SSS_{3}} e_{3}$
 does not depend on the choice of frames of the same orientation,
 and therefore is a well-defined section of the vector
 bundle $\wedge^2 TM\otimes \frak{su} (2)$.

 \end{thm}
 \pf  Assume that $\TT=\{e_{1}', e_{2}', e_{3}'\}$ is another
 frame and $(e_{1}', e_{2}', e_{3}' )=(e_{1}, e_{2}, e_{3})O$
 for some $O\in SO(3)$.
 It suffices to show that
 \begin{equation}
 (\pi_{\TT_{1}}, \pi_{\TT_{2}}, \pi_{\TT_{3}})=
 (\pi_{\SSS_{1}}, \pi_{\SSS_{2}}, \pi_{\SSS_{3}})O,
 \end{equation}
 or equivalently,
 \begin{equation}
 \label{eq:transformation}
 (\{f, g\}_{\TT_{1}}, \{f, g\}_{\TT_{2}}, \{f, g\}_{\TT_{3}})=
  (\{f, g\}_{\SSS_{1}}, \{f, g\}_{\SSS_{2}}, \{f, g\}_{\SSS_{3}})O
  \end{equation}
  for any $f, g\in C^{\infty}(M)$.

$M$ can be identified with $\tfrakg$ under both  $\Psi_{\SSS}$ and $\Psi_{\TT}$.
 If $(a,  b, c)$ and $(a', b', c')\in \tfrakg$ are, respectively,   the
coordinates of any point in $M$ under these  two  identifications,
they should be related by
$(a', b', c')=(a, b, c)O$. 
We assume that $O=\{(a_{ij})\}$. As an example, we  shall show below that
\begin{equation}
\{f, g\}_{\TT_{1}}=a_{11}\{f, g\}_{\SSS_{1}}+a_{21} \{f, g\}_{\SSS_{2}}
+a_{31} \{f, g\}_{\SSS_{3}},
\end{equation}
 for $f=<\xi, a>$ and $g=<\eta ,a>$. All the other cases can be 
proved similarly.

In this case, it is simple to see that
\be
\{f, g\}_{\TT_{1}}&=&(a_{11}, a_{12}, a_{13})\left( 
\begin{array}{lll}
<[\xi ,\eta], a'>&<[\xi ,\eta], b'>&<[\xi ,\eta], c'>\\
<[\xi ,\eta], b'>&-<[\xi ,\eta], a'>&A_{\xi, \eta}\\
<[\xi ,\eta], c'>&-A_{\xi, \eta}&-<[\xi ,\eta], a'>
\end{array}
\right)
\left( 
\begin{array}{c}
a_{11}\\a_{12}\\a_{13}
\end{array}
\right) \\
&=&(a_{11}, a_{12}, a_{13})\left(
\begin{array}{lll}
<[\xi ,\eta], a'>&<[\xi ,\eta], b'>&<[\xi ,\eta], c'>\\
<[\xi ,\eta], b'>&-<[\xi ,\eta], a'>&0\\
<[\xi ,\eta], c'>&0&-<[\xi ,\eta], a'>
\end{array}
\right)
\left(
\begin{array}{c}
a_{11}\\a_{12}\\a_{13}
\end{array}
\right)\\
&=&<[\xi, \eta], (a_{11}^2 -a_{12}^2-a_{13}^{2})a'+2a_{11}a_{12}b'+
2a_{12}a_{13}c' >\\
&=&<[\xi, \eta], (2a_{11}^2-1)a'+2a_{11}a_{12}b'+
2a_{12}a_{13}c' >\\
&=&<[\xi ,\eta], 2 a_{11}a-a'>.
\ee

On the other hand,
\be
&&a_{11}\{f, g\}_{\SSS_{1}}+a_{21} \{f, g\}_{\SSS_{2}}
+a_{31} \{f, g\}_{\SSS_{3}}\\
&=&<[\xi ,\eta], a_{11}a-a_{21}b-a_{31}c>\\
&=&<[\xi ,\eta], 2 a_{11}a-a'>.
\ee

This completes the  proof. \qed 

An immediate consequence  is the following:
\begin{thm}
 Assume that
 Equations (\ref{eq:conditionA}) and (\ref{eq:conditionB}) hold. Then the Poisson tensors
 $\pi_{\SSS}$,  $\pi_{\TT}$ corresponding to any  two  frames  commute.
 \end{thm}
 \pf For any $k_{1}, k_{2}, k_{3}\in \reals$ such
 that $k_{1}^2+ k_{2}^2+k_{3}^2=1$, it follows
 from  Theorem \ref{thm:global} that $k_{1}\pi_{\SSS_{1}}+ k_{2} \pi_{\SSS_{2}} 
 +k_{3}\pi_{\SSS_{3}}$
 is still  a Poisson tensor. Hence, $\pi_{\SSS_{1}}$, $\pi_{\SSS_{2}}$ 
 and $\pi_{\SSS_{3}}$ all commute.  Again according to  Theorem \ref{thm:global},
 $\pi_{\TT}$ is a linear combination of $\pi_{\SSS_{1}}$, $\pi_{\SSS_{2}}$ 
  and $\pi_{\SSS_{3}}$, and therefore commutes with $\pi_{\SSS} (=\pi_{\SSS_{1}})$. \qed
{\bf Remark} Under the assumption of Theorem   \ref{thm:conditionA},
the Poisson structure  $\pi_{\SSS} $ on $M$    is
$G$-invariant,  and $pr_{1}\smalcirc \Psi_{\SSS}: M \lon \frakg $,
the composition of $\Psi_{\SSS}$ with  $pr_{1}:
\tfrakg \lon \frakg$,  the projection onto its first factor,
is a $G$-equivariant momentum mapping. Similarly,
$\pi_{\SSS_{1}}$ and $\pi_{\SSS_{2}}$ are $G$-invariant
with equivariant momentum mappings $pr_{2} \smalcirc \Psi_{\SSS}$ 
and $pr_{3} \smalcirc \Psi_{\SSS}$, respectively.

The rest of the section  is devoted to   the investigation of 
the symplectic foliation of  $\pi_{\SSS}$.
For simplicity, whenever a frame $\SSS$ is fixed,
we shall omit the subscript $\SSS$ when
denoting the Poisson structure under the circumstance  without
confusion. I.e, we will  use $\pi_{i}$
to denote $\pi_{\SSS_{i}}$. 
By $X^{i}_{f}$, for $i=1, 2, 3$, we denote the vector 
field $\pi_{i}^{\#}(df)$ for $f\in C^{\infty}(M)$.
\begin{thm}
\label{thm:hypersymplectic}
The symplectic foliation of $\pi_{\SSS}$ is independent of the choice of 
frames $\SSS$,
and the induced family of
symplectic structures on each leaf is  hyper-symplectic
iff for any $(a, b, c)\in \tfrakg$ and $\xi \in \frakg$,
the following system of equations for $(u, v, w)$ has a 
solution:
\begin{equation}
\label{eq:hypersymplectic}
\left\{\begin{array}{lll}
 A^{\#} u+[v, c]-[w, b]&=&[\xi , a]\\
 -[u, c]+A^{\#} v+[w, a]&=&[\xi , b]\\
  +[u, b]-[v, a]+A^{\# }w&=&[\xi , c]\\
 -[u, a]-[v, b]-[w, c]&= &A^{\# } \xi,
 \end{array}
 \right.
 \end{equation}

\end{thm}
where for any $u\in \frakg$, $A^{\#}u$ is the
$\frakg$-valued function on $M$ obtained
by contracting with $u$, and
 $M$ is identified with $\tfrakg$ under  $\Psi_{\SSS}$.\\\\
{\bf Remark}\ \ The first three equations can be written in terms of
a single equation as 
$$X_{l^{3}_{u}}^2 +X_{l^{1}_{v}}^3+X_{l^{2}_{w}}^1=-\hat{\xi}. $$

If such an $A$ exists, we shall call the corresponding family
of Poisson structures   a  {\em hyper-Lie Poisson structure} and
their symplectic foliation {\em hypersymplectic foliation}.

We need several lemmas before we can prove this theorem.

The next lemma indicates that whether
System (\ref{eq:hypersymplectic})
is solvable is independent of the choice of frames. 
Therefore, the statement in Theorem \ref{thm:hypersymplectic}
is well justified.

\begin{lem}
\label{lem:indep}
For any fixed  $\xi \in \frakg$, $(u, v, w)$ is
a solution of System (\ref{eq:hypersymplectic})  for $(a, b, c)\in \tfrakg$
iff $(u', v', w')=(u, v, w)O$ is  a solution of the same
system for   $(a', b', c')$, where $O$ is any matrix
in $SO(3)$ and $(a' ,b', c')=(a ,b, c)O$.
\end{lem}

This can be  proved by a straightforward verification, and is  left
to the reader. The proof of the following
two lemmas is  also quite straightforward from definition.

\begin{lem}
\label{lem:pi1}
For the Poisson structure $\pi$ on $\tfrakg$, the hamiltonian vector fields 
of linear functions  are given by
\be
X_{l^{1}_{\xi}}&= &(-[\xi , a], -[\xi, b ], -[\xi, c]),\\
X_{l^{2}_{\xi}}&=&(-[\xi ,b], [\xi , a], A^\# \xi ),\\
X_{l^{3}_{\xi}}&=&(-[\xi ,c], -A^\# \xi , [\xi , a]),
\ee
\end{lem}
for any $\xi \in \frakg$, where
$(a, b,c)$ is any point in $\tfrakg$ and
 the tangent space at this point  is naturally  identified
with $\tfrakg$.
\begin{lem}
\label{lem:xij}
For any $\xi \in \frakg $ and $i, j=1, 2, 3$,
\be
X_{l_{\xi}^i }^j &=&-X_{l_{\xi}^j }^i , \ \ i\neq j;\\
X_{l_{\xi}^i }^i&=&-\hat{\xi}.
\ee
\end{lem}
{\bf Proof of Theorem  \ref{thm:hypersymplectic}. } Assume that System  
(\ref{eq:hypersymplectic}) has a solution.
We divide our proof into several steps.

(1). $\pi_{\SSS_{1}}^\# T^* M= \pi_{\SSS_{2}}^\# T^* M= \pi_{\SSS_{3}}^\# T^* M$.

From  Lemma \ref{lem:xij}, it follows that $X_{l_{\xi }^1}^2=-X_{l_{\xi }^2}^1$
and $X_{l_{\xi }^2}^2=-\hat{\xi}=X_{l_{\xi }^1}^1$. Hence, both
$X_{l_{\xi }^1}^2$ and $X_{l_{\xi }^2}^2$ are in $\pi_{\SSS_{1}}^\# T^* M$.
Also , it is easy to see that $X_{l_{\xi }^3}^2= (A^\# \xi , -[\xi , c], [\xi , b])$.
On the other hand,  for any $(u, v, w)\in \tfrakg \cong T^*_{(a ,b, c)}(\tfrakg)$,
 we have 
$$\pi_{\SSS_{1}}^\# (u, v, w)=(-[u, a]-[v ,b]-[w ,c],
-[u, b]+[v ,a]-A^\# w, -[u, c]+A^{\#}v +[w ,a])$$
 by Lemma  \ref{lem:pi1}.
It  is equal to $X_{l_{\xi }^3}^2$  if 
$(u, v, w)$ is a solution of  System  (\ref{eq:hypersymplectic}).
Hence, we have  $X_{l_{\xi }^3}^2\in \pi_{\SSS_{1}}^\# T^* M$. 
This shows that $\pi_{\SSS_{2}}^\# T^* M \subseteq \pi_{\SSS_{1}}^\# T^* M$.
Similarly $\pi_{\SSS_{1}}^\# T^* M\subseteq \pi_{\SSS_{2}}^\# T^* M$,
 so do the other relations as well.

(2). The symplectic foliation of $\pi_{\SSS}$
does not depend on the choice of frames.

Let $\TT$ be another frame, and $\pi_{\TT}$ its
corresponding Poisson structure.
According to Theorem \ref{thm:global}, $\pi_{\TT}$
can be expressed as a linear combination
of $\pi_{\SSS_{1}}$, $\pi_{\SSS_{2}}$ and $\pi_{\SSS_{3}}$. Therefore,
it follows that $\pi_{\TT}^\# T^*M\subseteq  \pi_{\SSS}^\#  T^*M$
from Step (1).
According to Lemma \ref{lem:indep} and
Step (1),
the symplectic foliations of $\pi_{\TT_{i}}, \ i=1, 2, 3$,  also
coincide. Hence,   exchanging $\SSS$ and $\TT$, we obtain  the other
inclusion: $\pi_{\SSS}^\# T^*M\subseteq  \pi_{\TT}^\#  T^*M$.

(3). Let $\omega_{1}$, $\omega_{2}$
and $\omega_{3}$ be the symplectic structures on a
hypersymplectic leaf  corresponding to
  $\pi_{\SSS_{1}}$, $\pi_{\SSS_{2}}$ and $\pi_{\SSS_{3}}$. Then  $\omega_{1}$, $\omega_{2}$, $\omega_{3}$
satisfy Equations (\ref{eq:compatible1}).

Below we only prove the following equation:  $\forall \ f, g\in C^{\infty}(M)$,
\begin{equation}
\label{eq:com}
\omega_{1}(X_{f}^2, X_{g}^2 )=\{f , g\}_{1}. 
\end{equation}
The other equations can  also be proved similarly.

In fact, it  suffices to show  this equation when 
 $f$ and $g$  are linear  functions.
If either $f$ or $g$ is $l_{\xi}^1 $ (for example $f=l_{\xi}^1$), we have 
\be
\mbox{LHS}&=&\omega_{1} (X_{l_{\xi}^1}^2, X_{g}^2)\ \mbox{ (by Lemma \ref{lem:xij})}\\
&=&-\omega_{1} (X_{l_{\xi}^2}^1, X_{g}^2)\\
&=&-<dl_{\xi}^2,  X_{g}^2>\\
&=&-\{g, l_{\xi}^2\}_{2} \\
&=&X^2_{l_{\xi}^2} (g) \ \mbox{ (by Lemma \ref{lem:xij})}\\
&=&-\hat{\xi} (g)\\
&=&\{l_{\xi}^1, g\}_{1}.
\ee

Similarly, Equation (\ref{eq:com}) holds if  either $f$ or $g$ is
 $l_{\xi}^2$. 

It remains to check Equation (\ref{eq:com}) for $f=l_{\xi}^3$ and $g=l_{\eta }^3$.
In this case, we know that  $X^2_{l_{\xi}^3}\per \omega_{1}=(u, v, w)$
according to the proof in Step (1), where $(u, v, w)$ is a solution of System 
(\ref{eq:hypersymplectic}) and  is considered as  a cotangent 
vector at $(a, b, c)\in \tfrakg$. We also know that
 $X_{l_{\eta }^3}^2= (A^\# \eta  , -[\eta  , c], [\eta  , b])$.
 Therefore,
 \be
\omega_{1} (X_{l_{\xi}^3}^2,  X_{l_{\eta }^3}^2)&=&(X_{l_{\xi}^3}^2\per \omega_{1}) ( X_{l_{\eta }^3}^2)\\
&=&<u, A^{\#} \eta >+<v, -[\eta , c]>+<w, [\eta , b]> \\
&=&<A^\# u+[v, c]-[w, b], \eta > \ \mbox{ (by Equation (\ref{eq:hypersymplectic}))}\\
&=&<[\xi , a], \eta >\\
&=&<[\eta ,\xi ], a>.
\ee

On the other hand,
  $\{l_{\xi}^3, l_{\eta }^3\}_{1}=l_{[\eta , \xi ]}^1=<[\eta ,\xi ], a>$.
  Hence, $\omega_{1} (X_{l_{\xi}^3}^2,  X_{l_{\eta }^3}^2)=\{l_{\xi}^3, l_{\eta }^3\}_{1}$.
  This completes the proof of Equation (\ref{eq:com}).

  Finally,   it is  quite transparent from the proof above that
  the assumption in  the statement of Theorem  
  \ref{thm:hypersymplectic} should  also be necessary. \qed

To end this section, we give the following result  which
reveals the connection between the $S^2  (\frakg )$-valued function $A$  and the
induced pseudo-metric on the  hypersymplectic leaves.
\begin{thm}
\label{thm:metric}
Under the assumptions as in Theorem \ref{thm:conditionA} and
Theorem \ref{thm:hypersymplectic}, the \peso-metric $g$ on each
hyper-symplectic leaf is $G$-invariant,  and for any $\xi \in \frakg$,
$g(\hat{\xi}, \hat{\xi})=-A_{\xi, \xi}$.
\end{thm}
\pf Since the Poisson structures $\pi_{\SSS_{1}}$, $\pi_{\SSS_{2}}$
and $\pi_{\SSS_{3}}$ are all $G$-invariant,
so are their induced symplectic structures on each hyper-symplectic leaf. 
Hence, the \peso-metric $g$ is $G$-invariant. According to
Equation (\ref{eq:g}),
\be
g(\hat{\xi}, \hat{\xi})&=&-\pi_{1}(\omeb_{2}\hat{\xi}, \omeb_{3}\hat{\xi})\\
&=&-\pi_{1}(-dl_{\xi}^2, -dl_{\xi}^3 )\\
&=&-\{l_{\xi}^2, l_{\xi}^3\}_{1}\\
&=&-A_{\xi, \xi}.
\ee
\qed

 We have seen  that the vector-valued
function $A$ plays a fundamental role in defining a hyper-Lie Poisson
structure.   The theorem above
 leads to some nondegenerate  criterion that   $A$ should  satisfy,
 i.e.,  $A_{\xi, \xi}=0$ iff $\hat{\xi }=0$.
The work of  Kronheimer \cite{Kronheimer:LMS}
\cite{Kronheimer:JDG1990} very much supports
the existence of  $A$ for compact semi-simple Lie algebras.
A satisfactory solution to this problem should 
provide us a symplectic approach, and therefore an intrinsic 
explanation, on the existence of hyper\kahler
structures on adjoint orbits. 
The work on this project is still in progress.
In the rest of the paper, instead we will consider the
case that
$\frakg=\frak{su}(2)$. This case can be handled 
relatively more
easily because of its special  character as  a three-dimensional Lie
algebra. However, we shall see that  certain nontrivial results, some of
which are already quite striking,  
 can be deduced  even  in  such a simple case.

\section{The case of $\frakg=\frak{su}(2)$}
From now on, we will work on the special case that
$\frakg=\frak{su}(2)$.  In this case, 
a  function $A$ can be   explicitly  constructed on an open
submanifold of $M$,
and  the corresponding
hyper-Lie Poisson structures are  studied 
under the general set-up  in  the previous section.
  
By $\Phi$, we denote   the function on $M$ defined by:
$$\Phi (a, b ,c)=<a, [b, c]>, \ \ \ \forall (a, b, c)\in \tfrakg. $$
Here again 
 $M$ is identified with $\tfrakg$ under some chosen frame. This
equation defines the   well known   Lie algebra $3$-cocycle corresponding to the 
Killing form. However, here we consider it as a function on $\tfrakg$ instead
of $\wedge^{3}\frakg$.
Clearly, $\Phi$ is independent of the choice of frames, provided
that they have the same orientation. Hence, $\Phi$ can still be
considered as a well-defined function on $M$.
Let $M_{o}$ be the open submanifold of $M$ consisting of all
points  where  $\Phi\neq 0$. In other words, $M_{o}$  consists
of triples $(a, b ,c)$ which are linearly independent.
Let $A: M_{o}\subset  \tfrakg \lon S^2 (\frakg )$ be the
map given by
\begin{equation}
\label{eq:A}
A(a, b,c)=\frac{1}{\Phi }([a, b]\otimes [a, b]+[b, c]\otimes [b ,c]
+[c, a]\otimes [c, a]), \ \ \forall (a, b, c)\in  M_{o} .
\end{equation}

It is not difficult to check that the rhs of Equation (\ref{eq:A})
is invariant under the natural action of $SO(3)$,
so $A$ can indeed be considered as 
  a  well-defined map from $M_{o}$ to $S^2 (\frakg )$. 

\begin{thm}

$A$ is $G$-equivariant and satisfies the condition in
Theorem  \ref{thm:conditionA}.
\end{thm}
\pf That $A$ is $G$-equivariant can be verified  directly.

We note that $A$ is uniquely characterized by the following relations:
\begin{equation}
\label{eq:character}
a\per A=[b, c], \ b\per A=[c, a], \mbox{ and }  c\per A= [a, b] 
\end{equation}
for any $(a, b, c)\in  M_{o}$.
Applying the vector field $X_{l_{\eta}^2}$ on both sides
of the equation $a\per A=[b, c]$, we obtain
$$X_{l_{\eta}^2}a \per A+a \per X_{l_{\eta}^2} A=
[X_{l_{\eta}^2}b, c]+[b, X_{l_{\eta}^2}c], $$
where both sides are  considered as  a $\frakg$-valued function on $\tfrakg$.
Using Lemma \ref{lem:pi1},    we have 
$$-[\eta, b]\per A+a \per  X_{l_{\eta}^2} A=
[[\eta , a], c]+[b, A^{\#} \eta ].$$
Hence, 
 $$a \per  X_{l_{\eta}^2} A=[[\eta , a], c]+[b, A^{\#} \eta ]+A^{\#}
[\eta , b]. $$
By contracting with $\xi\in \frakg$,  it follows that
\be
<a, \xi\per  X_{l_{\eta}^2} A> &=&<\xi, [[\eta , a], c]>+<\xi, [b, A^{\#} \eta ]>
+<\xi,  A^{\#} [\eta , b]>\\
&=&<\xi, [[\eta , a], c]>+<[\xi, b], A^{\#} \eta>+<\xi,  A^{\#} [\eta , b]>\\
&=&<\xi, [[\eta , a], c]>+A_{ \eta, [\xi, b]}+A_{\xi,  [\eta , b]}.
\ee
Thus
\be
<a, (\xi\per  X_{l_{\eta}^2} A-\eta\per  X_{l_{\xi}^2}A )>
& = &<\xi, [[\eta , a], c]>-<\eta , [[\xi , a], c]>\\
&=&<a, [[\xi, \eta ], c]>.
\ee
Using the other two identities in Equation (\ref{eq:character}), similarly 
we  deduce that 
\be
<b, (\xi\per  X_{l_{\eta}^2} A-\eta\per  X_{l_{\xi}^2}A) >&=&<b, [[\xi, \eta ], c]>, \\
<c, (\xi\per  X_{l_{\eta}^2} A-\eta\per  X_{l_{\xi}^2}A) >&=&0.
\ee
Since the Lie  algebra   $\frakg=\frak{su}(2)$ is three dimensional and $\Phi (a, b, c)\neq 0$,
$\{a, b, c\}$ constitutes a basis of  $\frak{su}(2)$ at any point
in $M_{o}$.  Equation  (\ref{eq:conditionA}) thus
 follows immediately, similarly for Equation (\ref{eq:conditionB}). \qed 

In fact,  $A$ also satisfies 
 the assumption as in Theorem \ref{thm:hypersymplectic}.

\begin{thm}
\label{thm:su2}
For the function $A$ defined by   Equation (\ref{eq:A}),
 System (\ref{eq:hypersymplectic}) always has a solution for any $\xi \in \frakg$
and $(a, b ,c)\in M_{o}$.  So $M_{o}$ has a hyper-Lie Poisson
structure. In  particular, its  symplectic leaves  are hyper-symplectic.
\end{thm}
\pf Fix any point $(a, b, c)\in M_{o}$. Since $\frak{su}(2)$ is
three-dimensional and System    (\ref{eq:hypersymplectic})
is linear with respect to $u, v, w$, it suffices to prove  this statement
 for  any three linearly independent $\xi \in \frakg =\frak{su}(2)$.
Therefore, it is sufficient to prove  this
for $\xi =a, b$, and $c$. For this purpose, one can check directly  that
$u=v=0, w=-b$ is a solution for $\xi =a$;
$u=-c, v=w=0$ is a solution for $\xi =b$;
and $u=0, v=-a, w=0$ is a solution for $\xi =c$. \qed

In fact, in this special case,  the corresponding hypersymplectic
 foliation can be described  quite 
   explicitly.  By  $X$, we denote the gradient vector field  of  
$\Phi$, where  $M$  is equipped with the standard metric
induced from the Killing form on $\frakg$. As  a frame is chosen and  $M$ is
identified with $\tfrakg$, 
 the vector field $X$ at any  point $(a, b, c)$ can be  written as 
\begin{equation}
X=([b, c], [c, a], [a, b]).
\end{equation}

Since both the standard metric on $M$ and the function $\Phi$ are
$G$-invariant, the gradient vector field $X$ is also $G$-invariant.
Therefore, it follows that 
$$[X, \hat{\xi}]=0, \ \ \ \forall \xi \in \frakg. $$

\begin{thm}
\label{thm:action}
The symplectic foliation of $\pi_{\SSS}$ on $M_{o}$
coincides
with the orbits of the  Lie algebra action of the direct product Lie algebra $\reals \times \frak{su} (2)$,
with $X$ and $\hat{\xi}, \xi \in \frakg$ being  its  generators.
\end{thm}
\pf The symplectic distribution at any point $(a, b, c)$ is
spanned by $\pi^\#_{\SSS}(T^*M)$. By identifying
$\frakg^*$ with $\frakg$, $T_{(a, b, c)}^*M$ is
identified with $\tfrakg$ as a vector space.
  To compute  the symplectic distribution,
it is sufficient to compute the image of a basis of
$\tfrakg$ under the map  $\pi^\#_{\SSS}$. Since $\frakg$ is three dimensional,
$\{a, b, c\}$ can be considered as a  basis of   $\frakg$. Hence, it suffices
to do the computation for its  corresponding   basis in $\tfrakg$.
Using  Lemma \ref{lem:pi1}, we have
$\pi^\#_{\SSS}(\xi, 0, 0)=-\hat{\xi}=-([\xi ,a], [\xi, b], [\xi, c])$.
It is also easy to see that  $\pi^\#_{\SSS}(0, a, 0)=\hat{b}$,
 $\pi^\#_{\SSS}(0, b, 0) =-\hat{a}$ and $\pi^\#_{\SSS}(0, c, 0)=X$;
 $\pi^\#_{\SSS }(0, 0, a)=\hat{c}$,
 $\pi^\#_{\SSS }(0, 0, b)=-X$ and $\pi^\#_{\SSS }(0, 0, c)=-\hat{a}$.
 This concludes  the proof of the theorem. \qed

 \begin{pro}
\label{pro:free}
  The  Lie algebra action defined    as in Theorem \ref{thm:action} 
 is locally free on $M_{o}$, so its  orbits are all 4-dimensional.
 \end{pro}
 \pf If not, there is $\xi \in \frakg $ and $k\in \reals$
 not all zero, such that
 $\hat{\xi }+kX=0$. If $k\neq 0$, it follows that $[b, c]=\frac{1}{k}[\xi ,a]$.
 Hence, $\Phi =<a, [b, c]>=<a, [\xi, a]>=0$. This contradicts to
 the definition of $M_{o}$. If $k=0$, we have $\hat{\xi }=
 -([\xi, a], [\xi ,b], [\xi , c])=0$. It thus follows that
 $\xi =0$ since $\{a, b, c\}$ is a basis of $\frakg$. \qed.

The coming  result gives us  a complete set of casimir functions for
 the  hypersymplectic foliation.

\begin{pro}
\label{pro:casimir}
 The following  functions  $<a, b>, <b, c>, <c, a>, 
 <a, a>-<b, b>$ and $<b ,b>-<c, c>$ form
 a complete set of casimirs for the Poisson structure $\pi_{\SSS}$
 on $M_{o}$.
 \end{pro}
 \pf It is simple to see that these functions are all $G$-invariant.
 To  show that they are casimirs, it suffices to show that
 they  are killed by the vector field $X$, which
 can be checked directly. It is also easy to see
 that these functions are all independent, so it 
 follows from dimension counting that this set of
 casimirs is complete. \qed

To end this section, we look at the induced metric on each hyper-symplectic
leaf.  The metric on the infinitesimal generators $\hat{\xi}$ of
the $G$-action is already given by Theorem  \ref{thm:metric}.
In order to describe the metric, we only need to know its evaluation
on the vector field $X$,  which is  the content of
the following:

 \begin{pro}
 \begin{equation}
 g(X, X)=-\Phi.
 \end{equation}
 \end{pro}
 \pf We already know that  $X_{l_{\xi }^3}^2=
(A^\# \xi , -[\xi , c], [\xi , b])$. Thus,  $\omeb_{2}X=(0, 0, a)$.
According to Lemma \ref{lem:pi1} and Lemma \ref{lem:xij}, we 
have   
 $\omeb_{3}X=(b, 0, 0)$. Here, in both equations,  the right hand sides
 are considered as elements in the cotangent
 space of $M$,  being identified with $\tfrakg$.  Therefore,
 $g(X, X)=-\pi_{1}(\omeb_{2}X, \omeb_{3}X ) =-\pi_{1}((0, 0, a), (b, 0, 0))
 =-<\pi_{1}^\# (0, 0, a),  (b, 0, 0)>
 =-<[c, a], b>=-\Phi $. \qed

The following consequence follows immediately from 
this result combining      with  Theorem  \ref{thm:metric}.
\begin{thm}
\label{thm:leaf}
When $\frakg=\frak{su}(2)$ and $A$ is defined by Equation
 (\ref{eq:A}), each hypersymplectic leaf of $M_{o}$
 is a 4-dimensional hyper\kahler manifold.
\end{thm} 

\section{Moduli spaces   of solutions to Nahm's equations and (co)adjoint orbits}
This is a  continuation of the last section.
When $M$ is identified with $\tfrakg $ by choosing a frame
 $\SSS$, the vector field $X$ is 
written as $X=([b, c], [c ,a ], [a ,b])$. Its flow  is
thus given by the following system of equations:
\begin{equation}
\label{eq:Nahm}
\left\{
\begin{array}{lll}
\dot{a}&=&[b, c]\\
\dot{b}&=&[c, a]\\
\dot{c}&=&[a, b].
\end{array}
\right.
\end{equation}
Such a system is called  Nahm's equations, which was studied
by Kronheimer \cite{Kronheimer:LMS} modelled on the  study of
general Nahm's equations made by Donaldson \cite{Donaldson:1984cmp}.
For any $x\in M$, we denote by $\phi_{t}(x)$ the flow generated by $X$
through the point  $x$. 
 We denote by  $S$  the set of all points
$x\in M$ such that the flow $\phi_{t}(x)$ converges as $t\lon -\infty$.
$S$ has  a natural foliation according to the limit points.
Note that all the limit points  are critical
points of $\Phi$. We  denote  by $C$ the critical
set of $\Phi$. 
The  $G$-action on $M$ leaves      $C$ invariant. For any  orbit
$\calo$, we let $S_{\calo}$ be the submanifold  of $S$ consisting
of all points in $S$ whose trajectory under the gradient vector  field
$X$ converges to a point in  $\calo$,  as $t\lon -\infty$, so that  in particular
$S_{0}$ be such a submanifold  corresponding
to the zero orbit.  It is clear that $S=\bigcup_{\calo}S_{\calo}$,
 where the sum is over all the $G$-orbits in  $C$.
Kronheimer proved, using gauge theory, for a general semi-simple Lie algebra 
that certain  $S_{\calo}$ are hyper\kahler  manifolds and  are diffeomorphic
to  adjoint orbits of $\frakg^{\complex}$ \cite{Kronheimer:LMS}
\cite{Kronheimer:JDG1990}. Below we will prove  this  result
for the special  case of $\frak{su}(2)$,  as a consequence
of the hyper-Lie Poisson structure on $S$. Our approach is
quite elementary,  and the family of  symplectic structures  
on each leaf  $S_{\calo}$ is rather transparent.

To start  with, let us introduce a  function $F$ on $M$  by
$$F(a, b, c)=<a, a>+<b, b>+<c, c>,$$
where $M$ is identified with $\tfrakg$ by choosing  a frame.
It is simple to see that $F$ is indeed a well-defined 
function on $M$.

\begin{lem}
\label{lem:lie-derivative}
\be
L_{X}\Phi&=&\|X\|^2, \ \ \mbox{and} \\
L_{X} F&=&6\Phi.
\ee
\end{lem}
\pf The first identity follows from a  general property of
a gradient flow.

As for the second one, we have

\be
L_{X}F&=&L_{X}(<a , a>+<b, b>+<c, c>)\\
&=&2<a, [b, c]>+2<b, [c ,a]>+2<c ,[a ,b]>\\
  &=& 6\Phi.
  \ee
  \qed

The following result is crucial for characterizing
the elements  in $S$.

\begin{pro}
\label{pro:converge}
\begin{enumerate}
\item If  $\phi_{t}(x)$ converges  as $t\lon -\infty$,
$x$ is either a critical point of $\Phi$ or $\Phi (x)>0$. 
In the latter  case,  we  in fact  have $\Phi (\phi_{t}(x) )>0$
for all $t$ whenever $\phi_{t}(x)$ is defined.
\item If  $\phi_{t}(x)$ does not converge as $t\lon -\infty$
or   is not defined for all  $t\leq  0$ 
(i.e., $X$ is incomplete in the $-\infty$ direction), 
 $\Phi$ cannot be always nonnegative along the flow.
\end{enumerate}
\end{pro}
\pf By Lemma \ref{lem:lie-derivative}, $\difft \Phi (\phi_{t}(x))=L_{X}\Phi =
\|X\|^2 \geq 0$.
Hence,  $\Phi (\phi_{t}(x))$  is an increasing function with respect to  $t$. 
If    $\phi_{t}(x)$ converges  as $t\lon -\infty$,  the
limit point must be  a critical point. 
However, the critical points of $\Phi$ are defined
by the system of equations:
\begin{equation}
[b, c]=[c, a]=[a, b]=0.
\end{equation}
So $\Phi $ vanishes at any critical point. This yields that $\Phi (x)\geq 0$.
If $\Phi (x)=0$, it follows that $\Phi (\phi_{t}(x) )=0$, for all
$t\leq  0$. By taking derivative, we have $\|X\|=0$.
Thus, $x$ is a critical point.
In the case that $\Phi (x)>0$, it is not difficult to see that $\Phi (\phi_{t}(x) )$ has to stay positive for 
 all $t$ whenever $\phi_{t}(x)$ is defined, otherwise $x$ will be  a critical
point according to  the same  argument above. 

If  $\phi_{t}(x)$ is not defined for all  $t\leq 0$,
it must be unbounded as $t$ approaches  to a finite number.
If  $\phi_{t}(x)$ is defined for all $t\leq 0$ but
does not converge as $t\lon -\infty$,  it must be unbounded as
$t$ is sufficiently negative  since $\Phi $ is a
real analytic function. In both cases, $F (\phi_{t}(x) )\lon \infty$,
as $t\lon -\lambda$ ($\lambda$ is either a positive number or
$\infty$). Assume that $\Phi$ is always nonnegative
along the flow.
It follows from Lemma  \ref{lem:lie-derivative}
that $\difft  F (\phi_{t}(x) )=6 \Phi  \geq 0$.
So  $ F (\phi_{t}(x) )\leq F(x)$ when $t\leq 0$, which 
is a contradiction.
This concludes  the proof. \qed

By  $M_{+}$, we denote the submanifold of $M$ consisting of
all points where $\Phi$ is positive.
The theorem above yields that $S-C$ is contained in $M_{+}$.
Moreover, the vector field $X$ is complete in $S-C$.
It is clear  that $S-C$ is invariant under the $G$-action, hence invariant
under the action of the product Lie algebra $\reals\times
\frak{su}(2)$ as defined in  Theorem \ref{thm:action}.
In other words,  $S-C$ is a  hyper-Poisson submanifold of $M_{o}$.
To extend this hyper-Poisson structure to entire  $S$, it  
suffices to extend $\pi_{\SSS}$ to the critical set $C$. For this,
one only needs to extend the vector valued function $A$ to
the  critical set $C$. Since $C$ is the limit set, a natural way to 
extend $A$ is to take its  limit along the flow $X$. This
is in fact how we derive the formula below.

When $M$ is identified with $\tfrakg$ under a chosen 
frame, a point 
$x_{0}=(a_{0} , b_{0}, c_{0})\in \tfrakg$ is a critical point
 iff $a_{0} , b_{0}, c_{0}$
are parallel. Hence,  for any critical point, we can
always choose a frame so that the critical point
is of the form $(a_{0}, 0, 0)$ for some $a_{0}\in \frakg$
under the identification:  $M\cong \tfrakg$ using  this frame.
Such a frame is called a standard frame. 
Clearly, the element $a_{0}$ is unique up to a sign.
We then define the function  $A$ on $C$ under a standard frame by:
\begin{equation}
\label{eq:extension}
A: C\lon S^{2} (\frakg ), \ \  (a_{0}, 0, 0)\lon \|a_{0}\| \sum I_{i}\otimes I_{i}
-\frac{1}{ \|a_{0}\|} (a_{0}\otimes a_{0}),
 \end{equation}
 where $I_{i}, i=1,2, 3, $ is an orthogonal basis for $\frakg\cong \frak{su}(2)$. We also let $A=0$ at $x=0$.
 $A$ is clearly well-defined on $C$.


To show that such an extension is  smooth, we need  to give
an alternate description of $S$,  which is much easier to deal with.

For any given  nontrivial  critical point $x_{0}$, 
let $\SSS_{x_{0}}$ be a  standard  frame 
such that  $x_{0}=(a_{0}, 0, 0)\in \tfrakg$
when  $M$ is identified with $\tfrakg$ under $\SSS_{x_{0}}$. 
We denote, by  $\Sigma_{x_{0}}$,
the subset of $M$ consisting of 
all the points $x=(a ,b ,c)\in \tfrakg ( \cong M$ under $\SSS_{x_{0}} )$,
satisfying the condition:
\begin{equation}
\label{eq:sigma}
<a, b>=<b, c> =<a, c>=0, <b, b>=<c,c>,  <a,a>-<b,b>=<a_{0}, a_{0}>, \mbox{ and }
\Phi \geq 0.
\end{equation}

It is easy to check that  this definition is  well-justified, i.e.,
does not depend on the choice of the  standard frame $\SSS_{x_{0}}$.
Also, we define $\Sigma_{0}$ as  the subspace  of $\tfrakg$ consisting
of all points $(a, b, c)$ such that
\begin{equation}
<a, b>=<b, c> =<a, c>=0,  <a,a>=<b,b>=<c,c>,  \mbox{ and }
\Phi \geq 0.
\end{equation}
It is clear that these relations are preserved under
the  transformation $(a ,b, c)\lon (a' ,b' ,c')=(a ,b, c)O$ for
any  $O\in SO (3)$. Therefore, $\Sigma_{0}$ can also be considered as
a subset of $M$.

\begin{lem}
For any nontrivial critical point $x_{0}$, $\Sigma_{x_{0}}=S_{G\cdot x_{0}}$.
\end{lem}
\pf $\Sigma_{x_{0}}$ is obviously a  closed submanifold of $M$.
It is clear that if $x\in \Sigma_{x_{0}}$, then $\phi_{t}(x)$ will stay
in $\Sigma_{x_{0}}$ for all the $t$ whenever the flow is defined,
since $X$ is tangent to $\Sigma_{x_{0}}$.
Since the intersection of $\Sigma_{x_{0}}$ with the hypersurface
$\Phi=0$ is contained in the critical set $C$, we conclude that
$\Phi (\phi_{t}(x))>0$ if $x$ is not a critical point in $\Sigma_{x_{0}}$.
Thus, according to Proposition \ref{pro:converge}, $\phi_{t}(x)$ exists
for all $t\leq 0$ and converges as $t\lon -\infty$.
Let us  assume that $y$ is  the limit point of $\phi_{t}(x)$. Then
$y$ is a critical point and $y\in \Sigma_{x_{0}}$.
 Assume that $y=(u ,v, w)$ under the  standard  frame
$\SSS_{x_{0}}$.  It is not difficult to see by using
Equation (\ref{eq:sigma}) that $v=w=0$ and $<u, u>=<a_{0}, a_{0}>$. The
latter implies that $y\in G\cdot x_{0}$. Hence
$\Sigma_{x_{0}}\subseteq S_{G\cdot x_{0}}$.

Conversely, assume that $x$ is any point in $S_{G\cdot x_{0}}$
and
$\phi_{t}(x)\lon y \in G\cdot x_{0}$ as $t\lon -\infty$.
Fix a standard  frame $\SSS_{x_{0}}$
 so that under it $x_{0}=(a_{0}, 0 ,0)$.
Then  under this  frame $y=(u, 0, 0)$ with  $u\in G\cdot a_{0}$.
Hence, $<u, u>=<a_{0}, a_{0}>$. Suppose that
$x=(a, b, c)$ under the frame $\SSS_{x_{0}}$.
From Proposition \ref{pro:casimir}, it thus follows that
$<a, b>=<u, 0>=0$, $<b, c>=<0, 0>=0$, $<a, c>=<u, 0>=0$,
$<b, b>-<c, c>=0$ and $<a, a>-<b, b>=<u, u>=<a_{0}, a_{0}>$.
That is, $x$ is in $\Sigma_{x_{0}}$. This completes the proof.  \qed

From this lemma, it follows that for any $x_{0}, y_{0}\in C$,
$\Sigma_{x_{0}}$ and $\Sigma_{y_{0}}$ are either disjoint, or equal.
In the latter case, $x_{0}$ and $y_{0}$ must lie in the same $G$-orbits
of  $C$.  For this reason, we shall  use $\Sigma_{\calo}$ to denote
the space $\Sigma_{x_{0}}$ for any $x_{0}\in \calo$.
The lemma above shows that $S_{\calo}=\Sigma_{\calo}$. 
In fact,  this is also valid when $\calo$ is the trivial orbit. 
\begin{lem}
$$S_{0}=\Sigma_{0}. $$
\end{lem}
\pf That $S_{0} \subseteq  \Sigma_{0}$ follows immediately from
the fact that  the functions: $<a, b>, <a, c>$,
$ <b, c>, <a, a>-<b, b>$, and $<b, b>-<c, c>$ are all preserved by
the vector field $X$. 

As for the other direction, let us assume that $(a, b, c)$ is any
nontrivial point in   $\Sigma_{0}$. By the definition of  $\Sigma_{0}$,
we can write $a=\lambda e_{1}, b=\lambda e_2$ and
$ c= \lambda e_{3}$, where $\lambda$ is a positive number and
$\{e_{1}, e_{2}, e_{3}\}$ satisfies the  standard
$\frak{su}(2)$-relation:
$[e_{1}, e_{2}]=e_{3}$, etc. Therefore, it is not difficult
to see that $a_{t}=-\frac{\lambda e_{1}}{\lambda t-1}, 
b_{t}=-\frac{\lambda e_{2}}{\lambda t-1}, c_{t}=-\frac{\lambda e_{3}}{\lambda t-1}$ is the flow  through the point $(a, b, c)$.
Obviously, it goes to zero as $t\lon -\infty$. That is,  $(a, b, c)\in S_{0}$.\qed

Combining the two lemmas above, we have
\begin{pro}
\label{pro:eq}
For any $G$-orbit $\calo$ in $C$, we have
$$S_{\calo}=\Sigma_{\calo}. $$
\end{pro}

Now we are ready to prove the  smoothness of the extension.
\begin{thm}
$\pi_{\SSS}$ is smooth when  restricted to
each  $S_{\calo}$ for
any nontrivial $G$-orbit $\calo$.
\end{thm}
\pf It suffices to show that the extension $A_{\xi, \xi}$ 
as defined by Equation   (\ref{eq:extension}) is smooth for any
$\xi \in \frakg$ in  a neighborhood of $\calo$ in  $S_{\calo}$.
Let $r$ be the norm of the  elements in  $\calo$.
Then by Proposition \ref{pro:eq}, under a standard  frame, points
$(a, b,c)$ in $S_{\calo}$ are characterized by
the equations:
$$ <a, b>=<b, c> =<a, c>=0, <b, b>=<c,c>,  <a,a>-<b,b>=r^2. $$
Therefore, when a  point $(a, b,c)$ is in $S_{\calo}$ but  not in
$\calo$, we can always write $a=\sqrt{\lambda^2 +r^2} e_{1},
b=\lambda e_{2}, c=\lambda e_{3}$, where $\lambda=\sqrt{<b ,b>}$
and $\{e_{1}, e_{2}, e_{3}\}$ is an orthonormal basis of
$\frak{su}(2)$ satisfying the standard relation. Write $\xi=\xi_{1}e_{1}+\xi_{2}e_{2}
+\xi_{3}e_{3}$. Thus a routine calculation 
yields that:
$$A_{\xi, \xi}(a,b,c)= \sqrt{\lambda^{2}+r^2 }(\xi_{2}^2 +\xi_{3}^2 )
+\frac{\lambda^2}{\sqrt{\lambda^{2}+r^2 }} \xi_{1}^2. $$
Substituting $\xi_{2}^2 +\xi_{3}^2 =  <\xi, \xi> -\xi_{1}^2$, we
have 
$$A_{\xi, \xi}(a,b,c)=\sqrt{\lambda^{2}+r^2 } <\xi, \xi> 
-\frac{r^2}{\sqrt{\lambda^{2}+r^2 }}\xi^2_{1}. $$
It is trivial to see that the  extension of $A$  as given by
Equation (\ref{eq:extension})
coincides with the equation above when $\lambda =0$.
Since both  $\lambda^2=<b, b>$ and $\xi_{1}=\frac{<\xi , a>}{<a, a>}$
are smooth functions on $S_{\calo}$, $A_{\xi, \xi}$ is
clearly smooth on $S_{\calo}$ as well. \qed

For any nontrivial orbit $\calo\subset C$, it is obvious
that $S_{\calo}$ is invariant under the $G$-action, as well as that of 
the additive group    $\reals$ generated by $X$.
Therefore, for any point $x\in S_{\calo}-\calo$, 
its  hypersymplectic leaf $\call_{x}$, defined as in the previous section,
is contained in $ S_{\calo}$. Since  $ \call_{x}$ is a 4-dimensional
manifold according to Proposition \ref{pro:free},  it  can be considered 
as an open neighborhood of $x$ in $ S_{\calo}$.  Clearly, 
$ S_{\calo}$ is a union of these leaves together with
their boundary $\calo$. As observed early, there is
a standard  frame $\SSS$ such that when $M$ is identified
with $\tfrakg$ under this frame, any point in $\calo$
is written as $x=(a, 0, 0)$ with $a\in \frakg$. In this
way, $\calo$ is naturally identified with a (co)adjoint orbit
of $\frakg$. Although there is an ambiguity for the choice of
the frame $\SSS$, such an adjoint orbit
is uniquely determined,  and the identification
is unique  up to a sign.
In the following,  we will fix any   such a  frame $\SSS$, and
denote the  Poisson structures  $\pi_{\SSS_{1}}, \pi_{\SSS_{2}}$
and $\pi_{\SSS_{3}}$ simply 
by $\pi_{1}$, $\pi_{2}$ and $\pi_{3}$, respectively, and the induced
symplectic structures on any leaf by $\omega_{1}, \omega_{2}$ 
and $\omega_{3}$ for simplicity.

\begin{thm}
For any nontrivial $G$-orbit $\calo\subset C$, the extended 
hyper-Poisson structure on $S$ induces a  hyper\kahler structure
on         $S_{\calo}$.
Furthermore,  if we choose a frame as in the observation above,
then $\calo$ is a
symplectic submanifold with respect to $\omega_{1}$. 
In fact,  it is a K\"ahler submanifold with respect to
 $I$ and  its K\"ahler metric coincides with
the one on $\calo$ when it is naturally identified
with a (co)adjoint orbit.
In the meantime,  $\calo$
 is   a lagrangian submanifold with respect to $\omega_{2}$ and
$\omega_{3}$.
\end{thm}
\pf It remains to consider points in $\calo$.  For any $x\in \calo$,
one can directly verify  that $\pi_{i}^{\#} T_{x}^*M=T_{x}S_{\calo}$,
for $i=1,2,3$, 
 by using a  local coordinate chart  under a chosen  standard frame $\SSS$.
Hence, the bivector fields $\pi_{i}, i=1,2,3,$ are all tangent to $S_{\calo}$
 and nondegenerate along $\calo$.
According to Theorem \ref{thm:su2}, the corresponding symplectic
structures $\omega_{i}$ are all compatible along $\calo$
by continuity, hence compatible  in entire  $S_{\calo}$.
Since the induced  metric is negative definite on $S_{\calo}-C$, it
is also   negative definite   along $\calo$.
Hence, the extended hyper-Poisson structure induces a
hyper\kahler structure on $S_{\calo}$.
The rest of the  conclusion can be verified directly, again by
 using local coordinates. \qed

In fact, $S_{\calo}$ is closely  related to  adjoint orbits
of $\frakg^{\complex}$. To see this,
 let $\SSS$ be any frame, and $pr_{12}\smalcirc \Psi_{\SSS}: M\lon \frakg^{\complex} $
   the  composition of the identification
$\Psi_{\SSS}: M\lon \tfrakg$ with the 
projection $pr_{12}$ as defined in the
proof of  Theorem \ref{thm:conditionA}.   According to  Proposition \ref{pro:poisson}, $pr_{12} \smalcirc \Psi_{\SSS}$ is 
a Poisson map with respect to $\pi_{\SSS}$. Hence, the image of
$S_{\calo}$ is a  Poisson submanifold of $\frakg^{\complex} $.
The following theorem indicates that in generic case
the image is in fact a single adjoint orbit.
A general result was proved by Kronheimer \cite{Kronheimer:LMS}
using gauge theory. However, our proof in this special
case  here is quite elementary and
only uses some well-known facts in symplectic geometry.
By $\calo_{12}$ we denote the adjoint orbit in 
$\frakg^{\complex} $ containing  the image $(pr_{12}  \smalcirc \Psi_{\SSS})(\calo )$.

\begin{thm}
\label{thm:adjointorbit}
If  $\calo_{12}$ is a regular orbit in $\frakg^{\complex}$,
 then $ ( pr_{12}    \smalcirc \Psi_{\SSS} ) (S_{\calo})=\calo_{12}$.
In fact, 
$pr_{12} \smalcirc \Psi_{\SSS} : S_{\calo}\lon \calo_{12}$ is a  symplectic diffeomorphism,
where  $ S_{\calo} $ is equipped with the 
 symplectic structure  induced from the Poisson 
structure $\pi_{\SSS}$ and $\calo_{12}$
is equipped with the  (co)adjoint orbit symplectic structure.
\end{thm}

We need some lemmas first.
\begin{lem}
\label{lem:bound}
 Under the map $pr_{12}\smalcirc \Psi_{\SSS}: S_{\calo}\lon \frakg^{\complex}$,
the inverse image of any bounded region
is bounded.
\end{lem}
\pf  According to Proposition \ref{pro:eq} and by the  definition
of $\Sigma_{\calo}$,
 there is a  standard frame $\TT$, under which any point
$(a', b', c')\in \tfrakg$  in $S_{\calo}$
is characterized by
 $$<a', b'>=<b', c'> =<a', c'>=0, <b', b'>=<c',c'>,  <a',a' >-<b',b'>=r^2, $$
where $r$ denotes  the norm  of elements in $\calo$.
Suppose that $O=(a_{ij})\in SO(3)$ is the
transformation matrix between the given frame $\SSS$ and 
this standard frame $\TT$.
That is, $(a, b, c)=(a', b', c')O$.
It is simple to see that $<a, a>=a_{11}^2 r^2+<b', b'>$,
$<b,b>=a_{12}^2  r^2+<b', b'>$, and $<c ,c>=a_{13}^2  r^2+<b', b'>$.
Therefore, if $<a, a>$ and $<b, b>$ are bounded, $<b', b'>$ has to be
bounded. This implies that $<c, c>$ is bounded as well. \qed

Two immediate consequences are the following:

\begin{cor}
\label{cor:proper}
The map $pr_{12}\smalcirc \Psi_{\SSS}: S_{\calo}\lon \frakg^{\complex}$ is a proper map.
\end{cor}

\begin{cor}
\label{cor:closed}
The image of $S_{\calo}$ under $pr_{12} \smalcirc \Psi_{\SSS}$ is closed.
\end{cor}
\pf Using the map $\Psi_{\SSS}$, we may identify
$M$ with $\tfrakg$. Suppose that under this identification
 $(a_{n}, b_{n}, c_{n})$ is a sequence
in  $S_{\calo}$ such that $a_{n}\lon a_{0}$ and
$b_{n}\lon b_{0}$. We need to show that $(a_{0}, b_{0})\in ( pr_{12} \smalcirc \Psi_{\SSS})(S_{\calo})$.
It follows from Lemma \ref{lem:bound} that $ c_{n}$ is bounded.
Therefore there exists  a convergent subsequence $c_{n_{k}}$,
whose  limit is denoted by  $c_{0}$. 
Then, $(a_{0}, b_{0}, c_{0})$ is in $S_{\calo}$ since
$S_{\calo}$ is closed. This concludes  the proof. \qed
{\bf Proof of Theorem \ref{thm:adjointorbit}} 
  Again, let us identify $M$ with $\tfrakg$ by
$\Psi_{\SSS}$. It is easy to see that at any point $(a, b,c )\in \tfrakg$,
 $T(pr_{12})X$ is the tangent vector,  at $a+ib$,   
 generated by the adjoint action $ad_{ia}$. Hence, the projection of
entire flow: $(pr_{12}\smalcirc \Psi_{\SSS}) ( \phi_{t}(x))$ lies in a single orbit
$\calc$ for all $t$. Since $\phi_{t}(x)$ converges to $\calo$ as
$t$ goes to $-\infty$,  then $(pr_{12}\smalcirc \Psi_{\SSS}) (\calo) 
\subseteq \bar{\calc}$. Therefore, $\calo_{12} \subseteq \bar{\calc}$.
Since $\calo_{12}$ is a regular orbit, it thus follows that
$\calo_{12}=\calc$.   This means  that the image  $(pr_{12} \smalcirc \Psi_{\SSS})(S_{\calo})$ is
contained in $\calo_{12}$ (this part of the argument is due to Kronheimer
\cite{Kronheimer:LMS}). 

On the other hand, according to Corollary \ref{cor:closed},
$(pr_{12} \smalcirc \Psi_{\SSS})(S_{\calo})$ is a closed Poisson submanifold in $\calo_{12}$.
Hence, it must be the entire  orbit $\calo_{12}$. That is, 
$ (pr_{12} \smalcirc \Psi_{\SSS}): S_{\calo}\lon \calo_{12}$ is onto. This map is
automatically  a submersion since it is a Poisson map.
By dimension counting, it must be a local diffeomorphism.
However, by Corollary \ref{cor:proper}, it is
a proper map. Therefore, it must be  a covering.
Since $\calo_{12}$ is simply connected, $ (pr_{12} \smalcirc \Psi_{\SSS}): S_{\calo}\lon \calo_{12}$
is  thus a diffeomorphism. \qed

Finally, we will show that $S_{0}-\{0\}$ is diffeomorphic 
to the nilpotent orbit of $\frak{sl}(2, \complex )$.

\begin{thm}
\label{thm:nilpotent}
$S_{0}-\{0\}$ is a hypersymplectic  leaf of $M_{o}$, and   therefore
 is a  hyper\kahler manifold. For  any frame $\SSS$, 
   $pr_{12} \smalcirc \Psi_{\SSS}$ is a symplectic diffeomorphism
between  $S_{0}-\{0\}$ and  the nilpotent orbit of $\frak{sl}(2, \complex )$,
where  $S_{0}-\{0\}$ is
equipped with the  symplectic structure  corresponding to $\pi_{\SSS}$
and the nilpotent orbit is equipped with the  standard 
coadjoint symplectic structure.
\end{thm}
 \pf  It is clear
 that $S_{0}-\{0\}$ is a union of hypersymplectic  leaves since
it is invariant under both  the $G$-action and the flow of $X$. 
 Each hypersymplectic leaf is 4-dimensional, and therefore must be
 open in $S_{0}-\{0\}$. Since $S_{0}-\{0\}$ itself
 is connected,  it must be a single hypersymplectic  leaf.
Thus, it is hyper\kahler according to   Theorem \ref{thm:leaf}.
By Proposition \ref{pro:eq}, under the identification
$\Psi_{\SSS}: M\lon \tfrakg$,  a point  $(a, b, c)\in \tfrakg$
is in $ S_{0}-\{0\}$ iff $a=\lambda e_{1}, b=\lambda e_{2}$ and
$c=\lambda e_{3}$ for some   standard orthonormal
basis $\{e_{1},  e_{2}, e_{3}\}$ of $\frak{su}(2)$, and $\lambda>0$.
Hence, its image under $pr_{12}$: $a+ib$,  is clearly
in the nilpotent orbit of $\frak{sl}(2, \complex )$.
A similar argument as in Corollary \ref{cor:closed}
shows that $ (pr_{12} \smalcirc \Psi_{\SSS}) (S_{0}-\{0\})$ is in fact closed.
Hence it has to be the whole nilpotent orbit since
$ pr_{12} \smalcirc \Psi_{\SSS}$ is a Poisson map.
Finally, it is quite obvious that $pr_{12} \smalcirc \Psi_{\SSS}$ is
injective on $S_{0}-\{0\}$. In fact, we always have
 $c=[a, b]/\sqrt{<a, a>}$.
This concludes  our proof. \qed
{\bf Remark} As we have seen, (co)adjoint orbits of $\frak{sl}(2, \complex )$ are 
related to the  points  in $M_{0}$ which have bounded
trajectories  (in the $-\infty$ direction) under the gradient vector  field
$X$ and are contained in $M_{+}$. However, according to Theorem  
\ref{thm:leaf},  there are other hypersymplectic leaves of
$M_{0}$ which are contained in $M_{-}=M_{0}-M_{+}$.
It would be interesting to explore  further the geometric structures for 
those leaves, and in particular the connection with the hyper\kahler
metrics on the  cotangent bundles of hermitian symmetric spaces of
noncompact type  studied  recently by  Biquard and Gauduchon \cite{BG}.


\begin{thebibliography}{99}
\bibitem{Atiyah}
M. F. Atiyah, Hyper-K\"ahler manifolds, 
{\em Collection: Complex Geometry and Analysis}, Springer-Verlag Lecture Notes in Mathematics,
{\bf 1422} (1990), 1--13.
\bibitem{AtiyahH:1988}
M. F. Atiyah and N. J. Hitchin, The geometry and dynamics of magnetic monopoles,
Princeton University Press, 1988.
\bibitem{Biquard}
O.~Biquard, Sur les \'equations de {N}ahm et la
structure de {P}oisson des alg\`ebres de {L}ie semi-simples complexes, 
{\em Math. Ann.}  {\bf 304},
(1996), 253-276.
\bibitem{BG}
O.~Biquard and P. Gauduchon, Hyperk\"ahler metrics on cotangent 
bundles of hermitian symmetric spaces, preprint.
\bibitem{Donaldson:1984cmp}
S.~K. Donaldson, Nahm's equations and the classification of monopoles,
{\em Commun. Math. Phys.} {\bf 96}, (1984), 387--407.
\bibitem{Hitchin}
N.~J. Hitchin, Hyper\kahler manifolds, 
{\em S\'eminaire Bourbaki} { 44\`eme ann\'ee, n. 748},
Ast\'erisque, {\bf 206}, (1992), 137-166.
\bibitem{HKLR}
N.~J. Hitchin, A. Karlhede, U. Linstrom and M. Rocek,
Hyper\kahler metrics and supersymmetry,
{\em Commun. Math. Phys.} {\bf 108}
(1987), 535-589.
\bibitem{K-SM:1990}
Y.~Kosmann-Schwarzbach and F.~Magri,
Poisson-{N}ijenhuis structures,
{\em Ann. Inst. H.~Poincar{\'e} Phys. Th{\'e}or.} {\bf 53}, (1990),
35--81.
\bibitem{Kovalev}
A.G..~Kovalev, Nahm's equations and complex adjoint orbits, preprint.
\bibitem{Kronheimer:LMS}
P.~B.~Kronheimer, A hyper-K\"ahlerian structure on coadjoint
      orbits of  a semisimple complex group,
{\em J. of LMS}, {\bf 42} (1990) 193--208.
\bibitem{Kronheimer:JDG1990}
P.~B.~Kronheimer, Instantons and the geometry of the nilpotent variety,
{\em J. Diff. Geom.} {\bf 32} (1990), 473--490.
\bibitem{MarsdenW}
J.~Marsden and A.~Weinstein, Reduction of symplectic manifolds with symmetry,
{\em Rep. Math. Phys.} {\bf 5} (1974), 121-129.
\bibitem{Penrose}
R.~Penrose, Nonlinear gravitons and curved twistor theory,
{\em Gen. Ralativ. Grav.} {\bf 7} (1976),
31-52.
\bibitem{Vergne}
M. Vergne, Instantons et correspondance de Kostant-Sekiguchi,
{\em C. R. Acad. Sci. Paris.} {\bf t 320} S\'erie I, (1995), 901--906.
 \bibitem{Weinstein:1983}
Weinstein, A. The local structure of {P}oisson manifolds,
{\em J. Diff. Geom.} {\bf 18} (1983), 523--557.
\end{thebibliography}
\end{document}